\documentclass[a4paper,10pt]{amsart}
\usepackage[utf8]{inputenc}

\usepackage[english]{babel}
\usepackage{amsmath,amssymb,verbatim,mathrsfs,latexsym,paralist}
\usepackage{graphicx}
\usepackage{epstopdf}
\usepackage{indentfirst}
\usepackage{amsthm}
\allowdisplaybreaks[4]

\newcommand{\Z}{\mathbb{Z}}
\newcommand{\C}{\mathbb{C}}

\newcommand{\mf}{\mathfrak}
\newcommand{\g}{\mf{g}}
\newcommand{\h}{\mf{h}}

\numberwithin{equation}{section}
\newtheorem{theorem}{Theorem}[section]
\newtheorem{proposition}[theorem]{Proposition}
\newtheorem{conj}[theorem]{Conjecture}
\newtheorem{lemma}[theorem]{Lemma}
\newtheorem{corollary}[theorem]{Corollary}

\theoremstyle{remark}
 
\theoremstyle{remark}
\newtheorem{rmk}[theorem]{Remark}
   
\title{On Reeder's Conjecture for Type $B$ and $C$ Lie algebras.}
\author{Sabino Di Trani}\address{Dipartimento di Matematica e Informatica, Universit\`a di Firenze.
}
\email{sabino.ditrani@unifi.it}
\address{\emph{The autor has been partially supported by GNSAGA - INDAM group.}}

\begin{document}

\maketitle

\begin{abstract}
In the paper we propose a proof of Reeder's Conjecture on the graded multiplicities of small representation in the exterior algebra $\Lambda \g$ for simple Lie algebras of type $B$ and $C$.
\end{abstract}

\section{Introduction}
Let $\g$ be a simple Lie algebra over $\C$. Fix a Cartan subalgebra $\h$ and let $\Phi$ be the associated root system with Weyl group $W$. We choose a  set of positive roots $\Phi^+$ associated to a simple system $\Delta$. Let $\rho$ be the corresponding Weyl vector and $\theta$ the highest root w.r.t. the standard partial order $\leq$ on $\Phi^+$. We will denote with $\Pi^+$ the set of dominant weights, $\omega_i$ will be the $i$-th fundamental weight.
 The adjoint action of $\g$ on itself induces a degree preserving action of $\g$ on $\Lambda \g$, the exterior algebra of $\g$. Then it is natural to study the irreducible components appearing in $\Lambda \g$ and their graded multiplicities. 
In \cite{Reeder} Reeder computes explicitly the multiplicities for a special class of irreducible representations $V_\lambda$, called \emph{small}, i.e. such that $\lambda$ is in the root lattice and $2 \alpha$  is not smaller than $\lambda$ in the dominant order for all positive roots $\alpha$.
 Furthermore, Reeder looks at the problem of determining the graded multiplicities of these representations, by firstly studying the simpler case of the adjoint representation and then trying to generalize the result to the other small modules.
He conjectured that the graded multiplicity of the small representation $V_\lambda$ can be computed reducing to a problem on finite group representations involving the $W$ representation on the zero weight space $V_{\lambda}^0$.
Let $H$ (resp $H^h$) be the space of $W$-harmonic polynomials (resp. of degree $h$) on $\h$, i.e. the polynomials annihilated by costant coefficients $W$-invariant differential operators with positive degree. 
\begin{conj}[Reeder]\label{RedC}
Consider the two polynomials 
 \begin{equation*}
 P(V_\lambda,\bigwedge \g, u)=\sum_{n\geq 0}\dim Hom_\g(V_\lambda,\bigwedge^n \g)u^n,
\end{equation*}
\begin{equation*}
 P_W(V^0_\lambda,H,x,y)=\sum_{n\geq 0}\dim Hom_W(V^0_\lambda,\bigwedge^k\h\otimes H^h)x^k y^h.
\end{equation*}
If $V_\lambda$ is a small representation of highest weight $\lambda$, then the following equality holds:
\begin{equation*}
 P(V_\lambda,\bigwedge \g, q)= P_W(V^0_\lambda,H,q,q^2)
\end{equation*}
\end{conj}
Curiously, this conjecture was implicitly proved for the case $A_n$ already before Reeder's paper was published, in the works of Stembridge \cite{St1} for the "Lie algebra" part  and by Kirillov, Pak and  Molchanov for the "Weyl group" part. 
Moreover, in \cite{Stembridge} many potentially useful tools for a case by case proof of the conjecture are introduced.  Of crucial importance for our work are some recursive relations for the coefficients $C_\lambda(t,s)$ in the characters expansion of the Macdonald kernels. These polynomial rational functions, specialized at $t=-q$ and $s=q^2$, give exactly the Poincar\'e polynomials for multiplicities of the representations $V_\lambda$ in $\Lambda \g$.
We propose a case by case proof of the conjecture for classical algebras of type $B$ and $C$ using the recursive relations of \cite{Stembridge} and closed formulae of \cite{GnS} for the Weyl group part.
More precisely, we will study the rational functions $C_\mu(q,t)$ using the so called "minuscule" and "quasi minuscule" recurrences proved by Stembridge in \cite{Stembridge}. These recurrences reduce the proof of the Reeder's Conjecture to solving an upper triangular system of linear equations with polynomial coefficients.
The first sections of the paper are dedicated to explain our tools and  the "Weyl group part" of the conjecture.
In Section 4 we prove Conjecture \ref{RedC} in the case of odd orthogonal algebras. Starting from Stembridge's minuscule recurrence  and using the combinatorics of weights and the action of the Weyl group we find some nice closed expressions for the coefficients of recursive relations.
Some further simplifications allows us to reduce to a two terms relation between the $C_\mu$. An inductive reasoning concludes the proof in this case. 
For the proof in case $C_n$, contained in the Section 5, we change our strategy and use the quasi minuscule recurrence.
Fixed a weight $\lambda$, the coefficients of the recurrence for $C_\lambda$  are described in terms of some suitable subsets in the orbit of $(\lambda, \theta)$ under the action of the Weyl group. We use the combinatorics of weights to express the associated coefficients in a recursive way. 
In the case of small weights of the form $\omega_{2k}$ the system of linear equations reduces easily to a two terms recursion that we solve using an inductive process.
For the weights $\omega_1 + \omega_{2k+1}$ the problem is more difficult: the zero weight representation $V_\lambda^0$ is not irreducible (except in the case $k=0$) and the combinatorics of the coefficients is more involved. 
We prove that the system of equation for $C_{\omega_1 + \omega_{2k+1}}$ can be reduced to a three terms relation involving $C_{\omega_{2k}}$ and $C_{\omega_{2(k+1)}}$. This allow us to reduce the conjecture to prove a univariate polynomials identity that we verified using SAGE.

\textbf{Aknowledgemets:} I am grateful to Professor Paolo Papi who suggested this problem for my PhD thesis and patiently supported me in the long and complicated process of review of this first article. Moreover I would like to thank Professors De Concini and Papi for sharing with me the sketch of their computations in case of odd orthogonal Lie algebras. Finally I am grateful to Professor Claudio Procesi for his advice on the organization of this paper.
\section{Stembridge's Recurrences}
Our main tools are the coefficients $C_\mu (q,t)$ in the character expansion of Macdonald kernels. Let  $\Delta(q,t)$ denote the Macdonald kernel and define $C_\mu (q,t)\in \C[q^{\pm 1},t^{\pm 1}]$ by the relation $\Delta(q,t)=\sum_{\mu\in  \Pi^+}  C_\mu (q,t) \chi(\mu)$. 
Extend the definition of $C_\mu (q,t)$ to any weight $\mu$ setting 
 \begin{equation}\label{RedClambda}
C_\mu(q,t) =
\left\{
	\begin{array}{lll}
	0 & \mbox{if } \mu + \rho \mbox{ is not regular}, \\
    (-1)^{l(\sigma)}C_\lambda(q,t)  & \mbox{if } \sigma(\mu + \rho) = \lambda+\rho \; , \lambda \in \Pi^+, \sigma \in W .\\
	\end{array}
\right.
\end{equation}
For short, we will say that, if there exists $\sigma$ such that $\sigma(\mu + \rho)=\lambda + \rho$, the weight $\mu$ is \emph{conjugated} to $\lambda$ and we will write $\mu + \rho \sim \lambda + \rho$.
In  \cite{Stembridge} Stembridge proves that the rational functions $C_\mu(q,t)$ satisfy some recurrences, reducing the problem of their explicit computation to solving a linear system of equations with coefficents in $ \C[q^{\pm 1},t^{\pm 1}]$.

We denote with $( \_\, ,\,\_ )$ the $W$-invariant positive-definite inner product on $\h^*$ induced by the Killing form and with $\alpha^\vee$ the coroot associated to $\alpha$. 
We recall that a weight (resp. coweight) $\omega$ is said to be \emph{minuscule} if $(\omega, \alpha^\vee ) \in \{0, \pm1\} $ (resp. $(\omega, \alpha)$) for all positive roots and \emph{quasi minuscule} if $(\omega, \alpha^\vee ) \in \{0, \pm 1 , \pm 2\} $ (resp. $(\omega, \alpha)$) for all positive roots.
Fix a dominant weight $\lambda$. If $\omega $ is a minuscule coweight, then the following relation holds (see \cite{Stembridge}, formula (5.14)):
 \begin{equation}\label{ricorsionemin}
  \sum_{i=1}^k C_{w_i \lambda}(q,t) \left(\sum_{\psi \in O_\omega} \left(t^{-(\rho, w_i \psi)}-q^{(\lambda, \omega)}t^{(\rho, w_i\psi)}\right)\right)=0.
 \end{equation} 
 Here $W_\lambda$ is the stabilizer of $\lambda$, $w_1, \dots , w_k$ are minimal coset representatives of $W/W_\lambda$ and $O_\omega$ is the orbit $W_\lambda \cdot \omega$. We are going to call this recursive relation the \emph{minuscule recurrence}.
 If $\g$ is not simply laced, the coroot $\theta^\vee$ is a quasi minuscule coweight. In such a case Stebridge proves the following recurrence: 
 \begin{equation}\label{ricorsioneqm}
\sum_{(\mu,\beta)}\sum_{i\geq 0}\left[ f_i^\beta(q,t)-q^{(\lambda,\omega)} f_i^\beta(q^{-1},t^{-1})\right]C_{\mu-i\beta}(q,t)=0.
\end{equation}
Here the pairs $(\mu,\beta)$ are elements of the set $\{(w\lambda,w \theta) | w \in W , \, w \theta \geq 0\}$  
and the rational functions $f_i^\beta$ are defined by:
\begin{equation}\label{fij}(1-tz)(1-qtz)\frac{((t^2z)^{(\rho,\beta^\vee)} -1)}{t^2z-1}=\sum_{i\geq 0} t^{(\rho,\beta^\vee)}f_i^\beta(q,t)z^i.\end{equation}
We will refer to (\ref{ricorsioneqm}) as the \emph{quasi-minuscule recurrence}.

 In both equations (\ref{ricorsionemin}) and (\ref{ricorsioneqm}) the polynomial rational functions $C_\mu$ are not necessary in their reduced form (i.e. the weight $\mu$ is not necessary dominant). The reduced form can be always achieved according to the Definition \ref{RedClambda}.
Considering only the reduced forms, Stembridge proves that the $C_\mu(q,t)$ appearing in (\ref{ricorsionemin}) and (\ref{ricorsioneqm}) are indexed only by weights $\mu$ smaller or equal to $\lambda$ in the dominant order. 
Our strategy is the following: we determine closed formulae for the polynomials of the "Weyl Groups" part of Reeder's Conjecture (explicitly computed in \cite{GnS}), and then we prove by induction that these closed formulae satisfy Stembridge's specialized recursive relations.
\section{Small Representations and Kirillov - Pak - Molchanov Formulae}
We recall that an irreducible representation $V_\lambda$ is \emph{small} if its weight is small.
The small representations for classical algebras and the structure of their zero weights space as $W$ representation are given in the simply laced cases by Reeder in \cite{R4}. In \cite{Kost} Kostant attributes the complete description of zero weight spaces to Chary and Pressler. 
We recall that the irreducible representations of hyperoctahedral group $S_n\ltimes \left(\Z/2\Z\right)^n$ (that by abuse of notation we are going to denote with $B_n$ as the associated root system) are encoded by pairs of partitions
$(\nu,\mu),\nu\vdash k,$ $\mu\vdash h, h+k=n$ and realized as $\pi_{\nu,\mu}={\mathrm{ Ind}}_{B_k\times B_h}^{B_n} \pi'_\nu\times \pi''_\mu$,  where, if $\pi_\tau$ is the irreducible $S_p$-module attached to 
$\tau\vdash p$, and $\varepsilon_q$ is the sign representation of  $\mathbb Z_2^q$,  we have $(\pi'_\nu)_{|S_k}=\pi_\nu\,, (\pi'_\nu)_{|\mathbb Z_2^k}=1_k, (\pi''_\nu)_{|S_h}=\pi_\mu, 
(\pi''_\mu)_{|\mathbb Z_2^h}=\varepsilon_h$.
The tables below sum up the relevant informations.
%
%

%
\begin{table}[h!]\label{smallweightsB}
  \begin{center}
   \caption{Zero weights space of small representation: Type B}
   \begin{tabular}{c|c}
    \textbf{Small Representation} & \textbf{Zero Weight Space } \\
    Highest weight & $(\alpha , \beta)$ description \\
    \hline
    $\omega_i$, $i<n$, $i=2k$ &  $((n-k),(k))$ \\
    $\omega_i$, $i<n$, $i=2k+1$ &  $((k),(n-k))$ \\
    $2 \omega_n$, $n=2k$ & $((k),(k))$\\
     $2 \omega_n$, $n=2k+1$ & $((k),(k+1))$\\
   \end{tabular}
\end{center}
\end{table}
\begin{table}[h!]\label{smallweightsC}
  \begin{center}
   \caption{Zero weights space of small representation: Type C}
   \begin{tabular}{c|c}
    \textbf{Small Representation} & \textbf{Zero Weight Space } \\
    Highest weight & $(\alpha , \beta)$ description \\
    \hline
    $2\omega_1$ &  $((n-1),(1))$ \\
    $\omega_{2i}$ &  $((n-i,i), \emptyset )$ \\
    $\omega_1 + \omega_{2i+1}$ , $i>0$ & $((n-i-1,i),(1)) \oplus ((n-i-1,i,1),\emptyset)$
   \end{tabular}
\end{center}
\end{table}
We remark that in the $B_n$ case the zero weight spaces are all irreducible; this does not happen for small representations $V_{\omega_1 + \omega_{2i+1}}$ in type $C$.

Now we will display the closed formulae that express explicitly the polynomials $P_W$  appearing in  the ``Weyl Group'' part of the Reeder Conjecture.
We will encode partitions $\lambda=(\lambda_1, \geq \lambda_2, \dots, \geq, \lambda_n)$ by Young diagrams, displayed in the English way.  $h(ij)$, $c(ij)$ are the hook length and the content of the box $(ij)$ respectively; $|\lambda|$ and $n(\lambda)$ will denote the quantities $\sum_{i=1}^n \lambda_i$ and $\sum_{i=1}^n(i-1)\lambda_i$.
\begin{theorem}[\cite{GnS}, Proposition 3.3]
 Let $m_1, \dots , m_n$ be the exponents of Weyl group $B_n$ and let $\pi_{\alpha, \beta}$ be the irreducible representation representation indexed by the pair of partitions $(\alpha, \beta)$.
 \begin{equation}\label{fromabrutatau}
P_W(\pi_{\alpha,\beta}; x,y)= x^{2n(\alpha)+2n(\beta)+|\beta|} \prod_{(i,j) \in \alpha} \frac{1+yx^{2c(ij)+1}}{1-x^{2h(ij)}} \prod_{(i,j) \in \beta} \frac{1+yx^{2c(ij)-1}}{1-x^{2h(ij)}}\prod_{i=1}^n \left(1-x^{m_i+1}\right),
\end{equation}
\end{theorem}
We can rearrange the formula (\ref{fromabrutatau}) obtaining:
\begin{equation}\label{formulatau}
P_W(\pi_{\alpha,\beta}; x,y)= \prod_{(i,j) \in \alpha} \frac{x^{2(i-1)}+yx^{2j-1}}{1-x^{2h(ij)}} \prod_{(i,j) \in \beta} \frac{x^{2i-1}+yx^{2(j-1)}}{1-x^{2h(ij)}} \prod_{i=1}^n \left(1-x^{m_i+1}\right).
 \end{equation}
\section{Odd Orthogonal Algebras}
We recall that the rooth system $B_n$ can be realized in an eucliden vector space with basis $\{e_1, \dots, e_n\}$, considering the set of vectors $\{\pm e_i \pm e_j\}_{i \neq j, i,j \leq n} \cup \{\pm e_j\}_{j \leq n}$. We fix a positive system of roots choosing the vectors of the form $\{e_i \pm e_j\}_{i\neq j, i,j \leq n} \cup \{e_j\}_{j \leq n}$; according to such a description, the fundamental weights are $\omega_i= e_1+ \dots + e_i$ if $i < n$ and $\omega_n= \frac{1}{2}(e_1+ \dots + e_n)$. Moreover the Weyl vector is $\rho = \frac{1}{2} \sum (2n-2j+1) e_j$.
Using the zero weight spaces description in Table 1, we obtain the following formulae for $P_W(V_{\omega_i}^0, q^2,q)$:
\begin{equation}\label{pippo}
 P_W(\pi_{((n-k)(k))}, q^2,q)=q^{2k-1}(q+1)\binom {n}{ k}_{q^4}\prod_{j=1}^{n-k}(1+q^{4j-1})
\prod_{l=1}^{k-1}(1+q^{4l-1})
\end{equation}
 if $i=2k$ and 
\begin{equation}\label{PWBn2}
 P_W(\pi_{((k)(n-k))}, q^2,q)=q^{2(n-k)-1}(q+1)\binom {n}{ k}_{q^4}\prod_{j=1}^{k}(1+q^{4j-1})
\prod_{l=1}^{n-k-1}(1+q^{4l-1})
\end{equation}
 if $i=2k+1$. Here the $q$-analogue $(n)_q$ is 
 $(n)_q=\frac{q^n-1}{q-1}=q^{n-1}+ \dots + 1$
and $q$-factorial and $q$-binomial are defined as 
$$(n)_q!=\prod_{i=1}^n(i)_q,\quad \binom{n}{m}_q=\frac{(n)_q!}{(m)_q!(n-m)_q!}.$$
Similar formulae hold for $V_{2 \omega_n}$ just distinguishing the cases $n=2k$ and $n=2k+1$.
We want prove inductively that these polynomials satisfy the minuscule recurrence (\ref{ricorsionemin}). In type $B$ a minuscule coweight is $\omega=e_1$. The stabilizer $W_{\omega_i}$ is isomorphic to $S_{i}\times B_{n-i}$ and $W_{\omega_i}(e_1)=\{e_1,\dots, e_i\}$. 
 We consider now the recurrence (\ref{ricorsionemin}) and make the evaluation $q \rightarrow -q$ and $t \rightarrow q^2$, obtaining

\begin{equation}\label{St2}
\sum_{i=1}^lC_{w_i\lambda}\sum_{j=1}^k (q^{-2(\rho,w_i\psi_j)}+q^{1+2(\rho,w_i\psi_j)})=0.
\end{equation}
Writing all the $C_\mu$ in their reduced form the recurrence can be rewritten as 
\begin{equation}\label{redSt2}
\sum_{\mu \leq \lambda}\Gamma_\mu^{n\lambda}(q)C_\mu(q) = 0,
\end{equation}
for some coefficients $\Gamma_\mu^{n\lambda}(q)$ (denoted as $\Gamma_h^{n,k}$ if $\lambda=\omega_k$ and $\mu=\omega_h$ ). Actually our purpose is to make more explicit the coefficients $\Gamma_h^{n,k}$.
Set 
\begin{align}\label{c}c_m&=\frac{1-q^{2m}}{1-q^2}q^{-2n+1}(1+q^{4n-2m+1}),\\
\label{B}b_m&=(q+1)q^{-2m+2}\frac{1-q^{4m-2}}{1-q^2}.\end{align}
Let $w$ be a minimal coset representative of $W/W_\lambda$. It can be checked easily that if  $\lambda=\omega_k$, then $C_{w \lambda}=C_{\lambda}$ if and only if $w=id$, obtaining that the coefficient of $\Gamma_\lambda$  in the reduced recurrence is exactly $c_k$.
Similar formulae hold for the coefficients $\Gamma^{nk}_\mu$.
We will start their explicit computation from the coefficient $\Gamma_0$.
Set $J(h,k,r)=\{((j_1, \dots, j_r)) \; | \; h<j_1<j_2-1< \dots < j_r-(r-1) \leq  k\}$. We set 
\begin{equation}
  \Gamma(h,k;r)=\sum_{\underline j\in J(h,k,r)}(q+1)\sum_{s=1}^r(q^{2(j_s-1)}+q^{-2(j_s-1)}).
 \end{equation}
We want prove now that $\Gamma_0^{n,k}$ is equal to
 \begin{equation}\label{Gamma0keven}
  \Gamma_0^{n,k}= (-1)^s\Gamma(1,n,s).
 \end{equation}
if $k=2s$, and to
 \begin{equation}\label{Gamma0kodd}
  \Gamma_0^{n,k}=(-1)^{s+1}\left[\binom{n-s-1}{s}b_1+\Gamma(2,n,s)\right].
 \end{equation}
 if $k=2s+1$.
 \begin{lemma}\label{zeroconj} Suppose $w \lambda$ is conjugated to $0$. Then:
 \begin{itemize}
  \item \emph{if $k$ is even},  $w \lambda$ has all the coordinates equal to zero except for $k/2$ pairs of consecutive coordinates of the form $(-1,1)$. Moreover $w \lambda $ is conjugated to $0$ by a permutation of sign $(-1)^{k/2}$.
  \item  \emph{if $k$ is odd} then $w \lambda$  has all the coordinates equal to zero, except for a choice of $(k-1)/2$ pairs of coordinates equal to $(-1,1)$ and for the last one that must be equal to $-1$. In this case $w \lambda $ is conjugated to $0$ by a permutation of sign equal to $(-1)^{(k-1)/2+1}$.
 \end{itemize}
\proof
  Observe that the coordinates of $w \lambda$ are all equal to $1$ or to $-1$. Then the coordinates of $w \lambda + \rho$ must be all positive except for the last one that can be equal to $-1/2$. So $w \lambda$ is conjugated to $0$ just by elements of the symmetric group $S_n$ and by sign change on the last coordinate. The lemma now follows by direct inspection.
 \endproof
 \end{lemma}
 As an immediate corollary we can compute the number of the weights giving a contribution to $\Gamma^{nk}_0$.
 \begin{corollary}
  Set $Conj^{n,k}_0=\{ w \omega_k \; | w \omega_k + \rho \sim \rho\}$. Then 
 \begin{equation*}\label{CardG0}
|Conj^{n,k}_0| =
\left\{
	\begin{array}{lll}
	\binom{n-\frac{k}{2}}{\frac{k}{2}} & \mbox{if } k \mbox{ is even}, \\
    \binom{n-\frac{k-1}{2}-1}{\frac{k-1}{2}} & \mbox{if } k  \mbox{ is odd.} \\
	\end{array}
\right.
\end{equation*}
\proof
By above lemma, the cardinality of $Conj^{nk}_0$ is equal to the number of choices of pairs of consecutive coordinates between $n$  or $n-1$ coordinates respectively when $k$ is even or odd. We can count the number of these possible choices in the following way: we can choose $k$ indices from the set $I=\{1 \dots n-k\}$ and expand our choices as a pair of consecutive indices.
\endproof
 \end{corollary}
Actually, by Lemma \ref{zeroconj}, if $k=2s$, we have $w \cdot (W_{\omega_k} \cdot e_1)= w \cdot \{e_1, \dots, e_k\} = \{-e_{j_1}, e_{j_1+1}, \dots , -e_{j_s}, e_{j_2+1}\} $. Set $v_{w}=(j_1, \dots , j_s)$ and $J(n,k)=\{v_{w} \, | \, w \in W/W_{\omega_k}\}$. The vector $v_{w}$ contributes to $\Gamma_0$ with a term equal to 
\begin{align*}
 \sum_{t=1}^s \left[q^{2(\rho, e_{j_t})}+q^{1-2(\rho, e_{j_t})}+q^{-2(\rho, e_{j_t+1})}+q^{1+2(\rho, e_{j_t+1})}\right] =
  \sum_{t=1}^s (q+1)\left[q^{2(n-j_t)}+q^{-2(n-j_t)}\right].
\end{align*}
Summing on the vectors $v \in J(n,k)$ we obtain 
\begin{equation}
 \Gamma_0^{n,k} = (-1)^s (q+1) \sum_{v \in J(n,k)} \sum_{t=1}^s \left(q^{2(n-j_t)}+q^{-2(n-j_t)}\right).
\end{equation}
The case $k=2s+1$ is very similar: we have only to observe that the last coordinate of $w \omega_k$ must be equal to $-1$ and choose the pairs of consecutive indices between $\{1, \dots n-1\}$. We must add  $b_1$ to the previous expression, coming from the contribution of $(w \omega_k)_n=-1$. We obtain
\begin{align*}
 \Gamma_0^{n,k} = (-1)^{s+1} \sum_{v \in J(n-1,2s)} \left[(q+1) \sum_{t=1}^s \left(q^{2(n-j_t)}+q^{-2(n-j_t)}\right)+  b_1\right]= \\
 (-1)^{s+1} (q+1) \sum_{v \in J(n-1,2s)} \sum_{t=1}^s \left(q^{2(n-j_t)}+q^{-2(n-j_t)}\right)+  \binom{n-k-1}{s} b_1.
\end{align*}
%
%
%
%
%
We obtain the closed formulae (\ref{Gamma0keven}) and (\ref{Gamma0kodd}) just observing that $J(n,k)=J(0,n-1,s)$ in the even case and $J(n,k)=J(0,n-2,s)$  in the odd one, and that if $v=(j_1, \dots j_s) \in J(0,n-2,s)$ (resp. $v \in J(0,n-1,s)$ ) then $v'=(n-j_s+1, \dots , n-j_1+1) \in J(2,n,s)$ (resp. $J(1,n,s)$).
Similar formulae hold for $\Gamma_i^{n.k}$: 
\begin{equation}
 \Gamma_i^{n.k} = \binom{n-i-s}{s}c_{i}+\Gamma(1,n-i;s)
 \end{equation}
if $k-i=2s$. Otherwise, if $k-i=2s+1$ 
\begin{equation}
 \Gamma_i^{n.k}=\binom{n-i-s-1}{s}(c_{i}+b_1)+\Gamma(2,n-i;s)
\end{equation}
and can be proved analogously to the $\Gamma_0^{n.k}$ case observing that 
$w \lambda$ is conjugated to $\omega_i$ if and only if 
 $w \lambda = (1, \dots , 1 , \gamma')$, $\gamma' \in Conj_0^{n-i,k-i}$.
We can then  write the recurrence in the following way :
\begin{align}\label{riscritta}C_mc_m&=\sum_{i=0}^{[\frac{m-1}{2}]}(-1)^{i}C_{m-2i-1}\Bigg (\binom{n-m+i}{i}(c_{m-2i-1}+b_1)+\Gamma(2,n-m+2i+1;i)\Bigg )\\\notag&+\sum_{i=1}^{[\frac{m}{2}]}(-1)^{i-1}C_{m-2i}\Bigg (\binom{n-m+i}{i}c_{m-2i}+\Gamma(1,n-m+2i;i)\Bigg ).\end{align}
The coefficients $\Gamma_i^{n,k}$ are now more explicit but less handy for a concrete computation. In Section 6 we prove the following proposition, which yields a more explicit and compact form for the coefficients.
\begin{proposition}
\begin{equation}\label{rec1}C_mc_m=\sum_{i=1}^{[\frac{m+1}{2}]}C_{m-2i+1}b_{i}+\sum_{i=1}^{[\frac{m}{2}]}C_{m-2i}b_{n-m+i+1}.
\end{equation}
\end{proposition}
The solution of minuscule recurrence can be found now by algebraic manipulations. 
%
%
%
\begin{theorem}\label{MP} Set $h=n-m/2$ if $m$ is even and $h=(m-1)/2$ is $m$ is odd. Then, if $h<n$, 
\begin{equation}\label{mp}C_m=\binom {n}{ h}_{q^4}\prod_{j=1}^h(1+q^{4j-1})\prod_{r=1}^{n-h-1}(1+q^{4r-1})(q^{2(n-h)}+q^{2(n-h)-1}).\end{equation}
\proof We start with the even case: $m=2k$; first we write \eqref{rec1}:
\begin{equation}\label{uno}C_mc_m=\sum_{i=1}^{k}C_{2(k-i)+1}b_{i}+\sum_{i=1}^{k}C_{2(k-i)}b_{n-2k+i+1}=\sum_{j=0}^{k-1}
C_{2j+1}b_{k-j}+\sum_{j=0}^{k-1}C_{2j}b_{n-k-j+1}.\end{equation}
Set,  for  $h>0$,
$$S_h=\begin{cases}\frac{1}{1-q^2}\binom {n}{ h}_{q^4}\prod_{j=1}^{n-h-1}(1+q^{4j-1})\prod_{r=1}^{h-1}(1+q^{4r-1})(q+1)^2\quad&\text{for $h>0$},\\
\frac{1}{1-q^2}\prod_{j=1}^{n-1}(1+q^{4j-1})(q+1)
\quad&\text{for $h=0$}.\end{cases}
$$
Recalling \eqref{B}, and using by induction (\ref{mp}) we have, for 
 $h>0$
\begin{align}\label{cane}C_{2h+1}b_{k-h}+C_{2h}b_{n-k-h+1}&=S_h\Big [(1+q^{4h-1})q^{2(n-k)+1}(1-q^{4(k-h)-2})\\\notag &+(1+q^{4(n-h)-1})q^{4h+2k-2n-1}(1-q^{4(n-k-h+1)-2})\Big ]\\\notag&=S_h(q^{4h}-q^{4(n-h)})(q^{2n-2k}+q^{2k-2n-1}).\end{align}
If $h=0$ we have 
\begin{align*}C_{1}b_{k}+C_{0}b_{n-k+1}&=S_0\Big[q^{2n-2k+1}(q+1)(1-q^{4k-2})+(1+q^{4n-1})q^{2k-2n}(1-q^{4n-4k+2})\Big ]\\
&=S_0(q-q^{4n+1})(q^{2n-2k}+q^{2k-2n-1}).\end{align*}
Hence, if we set $P_h=S_h(q^{4h}-q^{4(n-h)}),\ h>0, \ \ P_0=S_0(q-q^{4n+1})$, we can rewrite \eqref{uno} as 
\begin{equation}\label{due}C_mc_m=\sum_{j=0}^{k-1}
C_{2j+1}b_{k-j}+\sum_{j=0}^{k-1}C_{2j}b_{n-k-j+1}=(\sum_{h=0}^{k-1}P_h)(q^{2n-2k}+q^{2k-2n-1}).\end{equation}
Now, by induction, 
$$C_{2k-2}c_{2k-2}=(\sum_{h=0}^{k-2}P_h)(q^{2n-2k+2}+q^{2k-2n-3}).$$
Substituting into \eqref{due} we obtain
$$\sum_{j=0}^{k-1}
C_{2j+1}b_{k-j}+\sum_{j=0}^{k-1}C_{2j}b_{n-k-j+1}=(P_{k-1}+\frac{C_{2k-2}c_{2k-2}}{q^{2n-2k+2}+q^{2k-2n-3}})(q^{2n-2k}+q^{2k-2n-1}).$$
Now
$$c_{2k-2}=\frac{1-q^{4k-4}}{1-q^2}q^{-2n+1}(1+q^{4n-4k+5})=\frac{1-q^{4k-4}}{1-q^2}q^{4-2k}(q^{2n-2k+2}+q^{2k-2n-3}),$$
whence
$$\frac{c_{2k-2}}{q^{2n-2k+2}+q^{2k-2n-3}}=\frac{1-q^{4k-4}}{1-q^2}q^{4-2k},$$
and in turn
$$\sum_{j=0}^{k-1}
C_{2j+1}b_{k-j}+\sum_{j=0}^{k-1}C_{2j}b_{n-k-j+1}=(P_{k-1}+C_{2k-2}\frac{1-q^{4k-4}}{1-q^2}q^{4-2k})(q^{2n-2k}+q^{2k-2n-1}).$$
Now observe that 
$$c_{2k}=\frac{1-q^{4k}}{1-q^2}q^{-2n+1}(1+q^{4n-4k+1})=\frac{1-q^{4k}}{1-q^2}q^{2-2k}(q^{2k-2n-1}+q^{2n-2k})$$
so we are reduced to prove that 
\begin{equation}\label{tre}C_{2k}(1-q^{4k})q^{2-2k}=(1-q^2)P_{k-1}+C_{2k-2}(1-q^{4k-4})q^{4-2k}.\end{equation}
Divide both sides   of \eqref{tre} by 
$\frac{1}{1-q^2}\binom {n}{ k-1}_{q^4}\prod_{j=1}^{n-k}(1+q^{4j-1})\prod_{r=1}^{k-2}(1+q^{4r-1})(q+1)$. We get, for the r.h.s.
$$q (1 + q^{4 n - 4 k + 3}) (1 - q^{4 k - 4}) + (q + 1) (q^{4 k - 4} - 
    q^{4 n - 4 k + 4})=q+q^{4k-4}-q^{4n}-q^{5-4k+4n}=$$$$q(1+q^{4k-5})(1-q^{4n-4k+4}),$$
and the same for the l.h.s.
So the proof is completed in the even case. 
Consider now the odd case,  $m=2k+1$.
After a  suitable changes of variables we have
$$C_{2k+1}c_{2k+1}=\sum_{j=0}^{k}C_{2j}b_{k-j+1}+\sum_{j=0}^{k-1}C_{2j+1}b_{n-k-j}$$
Let us look at 
$$C_{2h+1}b_{n-k-h}+C_{2h}b_{k-h+1}$$
From Formula (\ref{cane}) we immediately deduce that for $h>0$,
\begin{align*}C_{2h+1}b_{n-k-h}+C_{2h}b_{k-h+1}&=S_h\Big [(1+q^{4h-1})q^{2k+1}(1-q^{4(n-k-h)-2})\\&+(1+q^{4(n-h)-1})q^{4h-2k-1}(1-q^{4(k-h+1)-2})\Big ]\\&=S_h(q^{4h}-q^{4(n-h)})(q^{2k}+q^{-2k-1})\end{align*}
For $h=0$,
\begin{align*}C_{1}b_{n-k}+C_{0}b_{k+1}&=S_0\Big[q^{2k+1}(q+1)(1-q^{4(n-k)-2})+(1+q^{4n-1})q^{2k}(1-q^{4k)+2})\Big ]\\
&=S_0(q-q^{4n+1})(q^{2k}+q^{-2k-1}).\end{align*}
Now we know that 
$$C_{2k}c_{2k}=(\sum_{h=0}^{k-1}P_h)(q^{2n-2k}+q^{2k-2n-1})$$
and also
$$\frac{c_{2k}}{q^{2n-2k}+q^{2k-2n-1}}=\frac{1-q^{4k}}{1-q^2}q^{2-2k},$$
Thus we deduce that
$$C_{2k+1}c_{2k+1}=C_{2k}b_1+C_{2k}\frac{1-q^{4k}}{1-q^2}q^{2-2k}(q^{2k}+q^{-2k-1})$$
Multiplying by $1-q^2$ we get
$$C_{2k+1}(1-q^{4k+2})q^{-2n+1}(1+q^{4n-4k-1})=C_{2k}((q+1)(1-q^2)+(1-q^{4k})q^{2-2k}(q^{2k}+q^{-2k-1}))=$$$$C_{2k}(q^{4k}+q)(q^{-4k}-q^2).$$
A simple computation shows that if we assume the formula for $C_{2k+1}$,
we get 
$$C_{2k+1}c_{2k+1}(1-q^2)=C_{2k}(q^{4k}+q)(q^{-4k}-q^2).$$
\endproof
\end{theorem}
Note that \eqref{mp} can be rewritten as 
\begin{align}
\label{prima}C_{2h+1}&=\binom {n}{ h}_{q^4}\prod_{j=1}^h(1+q^{4j-1})\prod_{r=1}^{n-h-1}(1+q^{4r-1})q^{2(n-h)-1}(q+1),\\
\label{seconda}C_{2k}&=\binom {n}{ k}_{q^4}\prod_{j=1}^{n-k}(1+q^{4j-1})\prod_{r=1}^{k-1}(1+q^{4r-1})q^{2k-1}(q+1).\end{align}
Making a direct comparison with the formulae (\ref{pippo}) and (\ref{PWBn2}) we obtain immediately that the Reeder's conjecture is verified for odd orthogonal algebras.
 \section{Symplectic Algebras}
For type $C_n$ the rooth system can be realized as $\{\pm e_i \pm e_j\}_{i\neq j, i,j \leq n} \cup \{\pm 2e_j\}_{j \leq n}$. A positive system is given by vectors $\{ e_i \pm e_j\}_{i\neq j, i,j \leq n} \cup \{ 2e_j\}_{j \leq n}$. The fundamental weights are then of the form $\omega_i= e_1+ \dots + e_i$, the coweight $e_1$ is quasi minuscule and the Weyl vector is $\rho =  \sum (n-j+1) e_j$. 
We have two different families of formulae for $P_W(V_{\lambda}^0, q^2,q)$, depending on the reducibility of representation $V_\lambda^0$.
We will  denote the polynomials $P_W(\pi_{((n-k,k); \emptyset)})(q^2,q)$ by $\textbf{C}_{k,n}$ and the polynomials $P_W(\pi_{((n-k-1,k); (1))}\oplus \pi_{(n-k-1,k,1), \emptyset})(q^2,q)$ by $\textbf{C}_{2|k,n}$. Similarly we will denote by $C_{k,n}(q,t)$ and $C_{2|k,n}(q,t)$ the rational functions appearing in the Stembridge's recurrences, omitting the variables if we are considering their specialized version. We have:
\begin{align*}\label{formulaPW}
\textbf{C}_{k,n}= q^{4k-1}(q+1)\binom{n}{k}_{q^4} \frac{\left(q^{4(n-2k+1)}-1\right)}{\left(q^{4(n-k+1)}-1\right)} \prod_{i=1}^{n-k}\left(q^{4i-1}+1\right) \prod_{i=1}^{k-1}\left(q^{4i-1}+1\right),
\end{align*}
\begin{align*}
\textbf{C}_{2|k,n} =  \frac{\prod_{i=1}^{n-k-1}(1+q^{4i-1}) \cdot \prod_{i=1}^k(q^4+q^{4i-1} ) \cdot \prod_{i=1}^n(q^{4i}-1) \cdot P_{k,n}(q^2,q) }{H(n-k,k+1) (1-q^2)},
\end{align*}
where  $H(n,k)= \prod_{(ij)} (1-q^{4h(ij)}) $ and we set
\[P_{k,n}(x,y)=(x+y)(x^{2(n-k+1)}-1)(x^{2(k+1)}-1)+(x^4+yx)(x^{2(n-k)}-1)(x^{2k}-1).\]
There are the following relations between the polynomials $\textbf{C}_{k,n}$ and $\textbf{C}_{2|k,n}$:
\begin{equation}\label{ricorsionePW}
\textbf{C}_{k+1,n}=\textbf{C}_{k,n} \frac{q^4(q^{4(n-k+1)}-1)(q^{4(n-2k-1)}-1)(q^{4k-1}+1)}{(q^{4(k+1)}-1)(q^{4(n-2k+1)}-1)(q^{4(n-k)-1}+1)}, 
\end{equation}
\begin{align}\label{T2kk+1}
 \textbf{C}_{2|k,n}=
 \notag \textbf{C}_{k+1,n}\frac{P_{kn}(q^2,q)}{q^4(1+q^{4k-1})(1-q^4)} \frac{(q^{4(n-2k)}-1)}{(q^{4(n-2k-1)}-1)(q^{4(n-k+1)}-1)}=\\\textbf{C}_{k+1,n}
 \frac{P_{kn}(q^2,q)}{(1+q^{4(n-k)-1})(1-q^4)} \frac{(q^{4(n-2k)}-1)}{(q^{4(n-2k+1)}-1)(q^{4(k+1)}-1)}.
\end{align}
%
%
%
%
%
%
\subsection{Weights of the form $\omega_{2k}$}
The quasi minuscule recurrence becomes: 
\begin{equation}\label{ricorsioneomegak}
\sum_{(\mu,\beta)}\sum_{i\geq 0}\left[ f_i^\beta(q,t)-q f_i^\beta(q^{-1},t^{-1})\right]C_{\mu-2i\beta}(q,t)=0.
\end{equation}
The pair $(\mu,\beta)$ is an element of the set $\{(w\lambda,w e_1) | w \in W , \, 2w e_1 \geq 0\}$  
and the $f_i^\beta$ are defined by equation (\ref{fij}).
Let us denote with $F_i^j$ the coefficient $f_i^{e_j}(q,t)-q f_i^{e_j}(q^{-1},t^{-1})$.
Observe that setting $A(i,j)=n-i-j+2$ we have $F_i^j=-F_{A(i,j)}^j$. More precisely if $A(i,j)=i$ we obtain $F_i^j=0$.
Set \[\Gamma_\nu^{i,j}=\{(w \lambda, \epsilon) \; w \in W , \, \epsilon \in \{\pm 1\} \; | \; w e_1 = e_j, \, C_{w(\lambda - 2ie_j)}=\epsilon C_\nu\}\]
 \begin{lemma}\label{simmetria1}
  There exists a bijection between $\Gamma_\nu^{i,j}$ and $\Gamma_\nu^{A(i,j),j}$ that sends a pair $(w \lambda, \epsilon)$ to a pair of the form $(w' \lambda , - \epsilon)$.
  \proof
  Let $(\mu_1, \dots , \mu_n)$ be the coordinates of $w \lambda$. By definition of $\Gamma_\nu^{i,j}$ we have $\mu_j= 1-2i$. The map
\begin{equation*}
\Psi(\mu) =
\left\{
	\begin{array}{lll}
        \mu_h  & \mbox{if } h\neq j, \\
		1-2(n-j-i+2) & \mbox{if } h=j,
	\end{array}
\right.
\end{equation*}
induces the desired bijection.
  \endproof
 \end{lemma}
We recall that, if $\lambda=\omega_{2k}$, considering the reduced form of rational functions $C_{\mu-2i\beta}(q,t)$, relation (\ref{ricorsioneomegak}) can be rewritten as 
$\sum_{i=0}^k \Lambda_i^{k,n} C_i(q,t) = 0$, for some coefficients $\Lambda_i^{k,n}$.
By Lemma \ref{simmetria1} we can reduce without loss of generality to consider only the pairs $(i,j)$ such that $2i-1 < n-j+1$, i.e. all the coordinates of $w (\omega_{2k}-2ie_1)+\rho$ are positive. In this case $w (\omega_{2k}-2ie_1)+\rho$ can be rearranged in the form $\omega_h + \rho$ just using elements of the symmetric group $S_n$ acting on the coordinates.
As an immediate consequence, the only contributions to the coefficient $\Lambda_k^{kn}$ come from the case $i=0$ and a pair $\{(\mu, e_j)\}$ appearing in the quasi minuscule recurrence (\ref{ricorsioneomegak}) and giving contribution to $\Lambda_k^{k, \,n}$ must be of the form $(\omega_{2k},e_j)$ with $1 \leq j\leq 2k $. The following closed formula follows: 
\begin{equation}\label{coefficientek}
\Lambda_k^{k, \,n}=\sum_{i=1}^{2k}F_0^{i, \, n}=\sum_{i=1}^{2k} \left( \frac{1}{t^{n-i+1}}-qt^{n-i+1} \right)=\frac{(t^{2k-1}-qt^{2n})(t^{2k}-1)}{t^{n+2k-1}(t-1)}
\end{equation}
It is more difficult to obtain closed formulae for the generic coefficient $\Lambda_{h}^{k,n}$, however some nice recurrences hold. For the weights conjugated to 0 a result similar to the one proved for $B_n$ holds:
\begin{lemma}\label{zerocontribution}
Set $\lambda=\omega_{2k}$ and let $w\in W$ be such that $w \lambda$ is conjugated to $0$. 
\begin{enumerate}
 \item The $2k$ non zero coordinates of $w\lambda$ are pair of consecutive coordinates $((w\lambda)_{(j)},(w\lambda)_{(j)+1})$  of the form $(-1,1)$.
 \item There exists a permutation $\sigma \in S_n$ of length $l(\sigma)=k$ such that $\sigma (w \lambda + \rho)= \rho$.
\end{enumerate}
Moreover the number of weights of the form $w \lambda$ conjugated to $0$ is $\binom{n-k}{k}$.
\end{lemma}
Set $\Omega_{h }^{k, \; n}=\{ w \omega_{2k} \;| \;  w \omega_{2k} + \rho \sim \omega_{2h}+\rho \}$ and $O_{h }^{k, \; n}=\{ w (\omega_{2k} - 2je_1) \;| \;  w (\omega_{2k} + 2je_1) + \rho \sim \omega_{2h}+\rho \}$.
\begin{rmk}\label{omegah+zero}
A direct inspection shows that the weights giving non zero contribution to the coefficient $\Lambda_{h}^{k,n}$ for $h>0$ are of the form $e_1+\dots+e_{2h}+\mu$, where $\mu$ has the first $2h$ coordinates equal to 0. Considering the immersion of $C_{n-2h} \rightarrow C_n$ induced by Dynkin Diagrams, this means that $\mu$ can be contracted to a weight in $\Omega_{0}^{k-h,n-2h}$. By abuse of notation, we will denote this contraction process writing $\mu \in \Omega_{0}^{k-h,n-2h}$.
\end{rmk}
Actually, the contribution to $\Lambda_{h}^{k,n}$ depends by $w(1)$:
\begin{itemize}
 \item \emph{Case 1: $w(1)=j$ with  $j \leq 2h$.} In this case $i$ is forced to be $0$ because, as just observed, $\lambda_j $ must be equal to 1. By the previous remark, the rational function $F_0^j$ appears as many times as the cardinality of $\Omega_{0}^{k-h,n-2h}$ with a sign equal to $(-1)^{k-h}$ by Lemma \ref{zerocontribution}. Consequently the contribution to $\Lambda_h^{k, \, n}$ is equal to \[(-1)^{k-h} F_0^j \binom{n-h-k}{k-h}.\]
 \item \emph{Case 2: $w(1)=j$ with $j > 2h$.} By Remark \ref{omegah+zero}, the contribution in this second case is the same as considering the one given by a weight in $O_{0}^{k-h,n-2h}$. It is then equal to $\Lambda_0^{k-h, \, n-2h}$.
\end{itemize}
Summing up we obtain the relation
\begin{equation}\label{Lh}\Lambda_h^{k, \, n}= \Lambda_0^{k-h, \, n-2h} + (-1)^{k-h} \sum_{j=1}^{2h}F_0^j \binom{n-h-k}{k-h} =\Lambda_0^{k-h, \, n-2h} + (-1)^{k-h} \Lambda_h^{h, \, n} \binom{n-h-k}{k-h}   \end{equation}
For $\Lambda_0^{k,n}$ the argument is just a bit different. If we set 
$O_{0, \, i}^{k, \; n}=\{ w (\lambda - 2je_1) \;| \;  w (\lambda - 2je_1) + \rho \sim \rho, \, j=i \}$
and denote by $\Lambda_{0,i}^{k,n}$ the relative contribution to $\Lambda_0^{k,n}$. Consider $\mu \in O_{0, \, i}^{k, \, n} $, the only relevant cases are:
\begin{itemize}
 \item \emph{Case 1: $\mu_1=0$.} In this case $\mu$ can be contracted to a weight $\mu' \in O_{0, \, i}^{k, \; n-1}$
  \item \emph{Case 2: $i>0$ and $\mu_1=1-2i$.} 
 This forces $\mu$ to be of the form $(1-2i,1, \dots, 1, \mu')$ where $\mu'\in \Omega_{0}^{k-i,n-2i}$. 
  \item \emph{Case 3: $\mu_1=-1$ and $w(1)\neq 1,2$.} This happens only if $\mu_2=1$ and $\mu=(-1, 1 , \mu')$, where $\mu'$ is a weight in $O_{0, \, i}^{k-1, \; n-2} $. 
  \item \emph{Case 4: $\mu_1=-1$ and $w(1)=2$.} The only possibility is that $i=0$ and $\mu=(-1, 1 , \mu')$, with $\mu' \in\Omega_{0, }^{k-1, \; n-2} $. 
\end{itemize}
%
%
%
%
The previous analysis leads us to the relation
 \begin{equation}\label{L0}\Lambda_{0}^{k, \, n}= \Lambda_{0}^{k, \, n-1} - \Lambda_{0}^{k-1, \, n-2}+ (-1)^k\binom{n-k-1}{k-1}F_0^{2, n}+ \sum_{i=1}^{k}(-1)^{k-i+1} F_i^{1, n} \binom{n-i-k}{k-i}.  \end{equation}
 The above relation between coefficients allow us to reduce the triangular system given by Stembridge's relations as described in the following proposition, that we are going to prove in Section 7.
 \begin{proposition}\label{redrecursion}
 Let $R_i$ be the recurrence for $C_{i,n}(q,t)$ written in reduced form. Then there exist a family of integers $\{A_i^{k, n}\}_{i \leq k}$ such that 
 \begin{equation}\label{Riduzionebella}\sum_{i=1}^k A_i^{k,n}R_i = \Lambda_k^{k, n} C_{k,n}(q, t)+ \Lambda_0^{1, n-2k+2} \left(C_{k-1, n}(q,t)+ \dots + C_{0, n}(q,t)\right)\end{equation}
\end{proposition}
Using \eqref{Riduzionebella}, some explicit formulae for $C_{k+1,n}(q,t)$ in function of $C_{k,n}(q,t)$ can be obtained:
           \begin{theorem}\label{ReederC2k}
           With notation as in Proposition \ref{redrecursion}, the following formula holds:
\begin{equation}\label{odiolatex}
           C_{k+1, n}(q,t)= \frac{(t^{2(n-2k-1)}-1)(t^{2(n-k+1)}-1)(1-qt^{2k-1})t^2}{(t^{2(n-2k+1)}-1)(t^{2(k+1)}-1)(1-qt^{2(n-k)-1})} C_{k, n}(q,t).
            \end{equation}
          \proof
          By Proposition \ref{redrecursion} we know that 
           \[ \Lambda_{k+1}^{k+1, n} C_{k+1,n}(q, t)=- \Lambda_0^{1, n-2k} \left(C_{k, n}(q,t)+ \dots + C_{0, n}(q,t)\right)\]
 Multiplying by $ \Lambda_0^{1, n-2k+2}$ and using the
inductive hypothesis we obtain 
          \begin{align*}
           \Lambda_0^{1, n-2k+2}\Lambda_{k+1}^{k+1, n} C_{k+1,n}(q, t)=\\
           - \Lambda_0^{1, n-2k}\Lambda_0^{1, n-2k+2} C_{k, n}(q,t) - \Lambda_0^{1, n-2k} \left[\Lambda_0^{1, n-2k+2} \left(C_{k-1, n}(q,t)+ \dots + C_{0, n}(q,t)\right)\right]= \\
            - \Lambda_0^{1, n-2k}\Lambda_0^{1, n-2k+2} C_{k, n}(q,t) + \Lambda_0^{1, n-2k}\Lambda_k^{k, n} C_{k, n}(q,t) = \\
            - \Lambda_0^{1, n-2k} \left(\Lambda_0^{1, n-2k+2} - \Lambda_k^{k, n}\right)C_{k, n}(q,t),
          \end{align*}
           and then 
           \begin{align*}\label{relCk}
            C_{k+1,n}(q, t)=  - \frac{ \Lambda_0^{1, n-2k} \left(\Lambda_0^{1, n-2k+2} - \Lambda_k^{k, n}\right)}{ \Lambda_{k+1}^{k+1, n}\Lambda_0^{1, n-2k+2}}C_{k, n}(q,t).
           \end{align*}
\begin{lemma}
      \begin{equation}\label{blabla}
     \Lambda_0^{1, n}=-\frac{(t-q)(t^{2n-2}-1)}{t^{n-1}(t-1)}.
      \end{equation}
     \proof
     By Lemma \ref{zerocontribution}, in this case the weights conjugated to 0 are of the form $(0 \dots 0, -1, 1, 0 \dots , 0)$. This implies immediately that  $i \leq 1$. If $i=0$ then $w(1) \in \{2, \dots , n\}$ and the contribution to $\Lambda_0^{1, \, n}$ is equal to 
     $-\sum_{j=2}^n F_0^{jn}$.
If $i=1$ then $w(1) \in \{1, \dots , n-1\}$ and the contribution is equal to $-\sum_{j=1}^{n-1} F_1^{jn} $.
Recalling that $F_0^{j,n}=F_0^{2, n-j+2}$ and $F_1^{j-1,n}=F_1^{1, n-j+2}$ we obtain 
\[\Lambda_0^{1, \, n}=  -\sum_{j=2}^n F_0^{j, n} + F_1^{j-1, n} = -\sum_{j=2}^{n-1} \left( F_0^{2, j} + F_1^{1, j}\right)=-\frac{(t-q)(t^{2n-2}-1)}{t^{n-1}(t-1)}\]
     \end{lemma}
 Replacing the coefficients with their closed formulae we obtain the statement of Theorem \ref{ReederC2k}.
          \endproof
          \end{theorem}         
It is now easy to verify that Reeder's Conjecture holds for the weights of the form $\omega_{2k}$ just
specializing the expression (\ref{odiolatex}) at $q \rightarrow -q$   and $t \rightarrow q^2 $  and checking that we obtain exactly the formula (\ref{ricorsionePW}).
\begin{rmk}
  We remark that this computation allow us to obtain explicit closed formulae for the non specialized rational functions $C_\mu(q,t)$ when $\mu$ is of the form $\omega_{2k}$. 
\end{rmk}

\subsection{Weights of the form $\omega_1 + \omega_{2k+1}$} We remark that the case $k=0$ is well known and has been proved for example in \cite{DCPP} and in \cite{Stembridge}. We will proceed by induction.
The recurrence in this case becomes: 
\begin{equation}
\sum_{(\mu,\beta)}\sum_{i\geq 0}\left[ f_i^\beta(q,t)-q^2 f_i^\beta(q^{-1},t^{-1})\right]C_{\mu-2i\beta}(q,t)=0.
\end{equation}
Similarly to the previous case we will denote $F_i^j=f_i^{e_j}(q,t)-q^2 f_i^{e_j}(q^{-1},t^{-1})$ and $\Gamma_\nu^{i,j}=\{(w \lambda, \epsilon) \; w \in W , \, \epsilon \in \{\pm 1\} \; | \; w e_1 = e_j, \, C_{w(\lambda - 2ie_j)}=\epsilon C_\nu\}.$
\begin{lemma}\label{simmetria11k}
 Set $B(i,j)=n-i-j+3$. We have a bijection between $\Gamma_\nu^{i,j}$ and $\Gamma_\nu^{B(i,j),j}$.
This bijection sends a pair $(w \lambda, \epsilon)$ in a pair of the form $(w' \lambda , - \epsilon)$.
  \proof
Analogously to the previous case, the map defined as
\begin{equation*}
\Psi(\mu) =
\left\{
	\begin{array}{lll}
        \mu_h  & \mbox{if } h\neq j, \\
		2-2(n-j-i+3) & \mbox{if } h=j,
	\end{array}
\right.
\end{equation*}
gives the desired bijection.
  \endproof
 \end{lemma}
We will denote by $\Psi_i^{j,n}$ the expression $ F_i^{j,n}-F_{B(i,j)}^{j,n}$.
The recurrence (\ref{ricorsioneqm}) can be written in reduced form as:
\begin{equation*}
  R_k  \; :  \; \sum_{i=0}^k C_{2|i}(q,t)\Lambda^{2|k,n}_{2|i}+\sum_{i=0}^{k+1/k} C_{i}(q,t)\Lambda_i^{2|k, \, n}=0
 \end{equation*}
 for some coefficients $\Lambda^{2|k,n}_{2|i}$ and $\Lambda^{2|k,n}_{i}$ that we are going now to analyze more closely. Moreover the index of the second sum goes from 0 to $k$ if $n=2k+1$ and to $k+1$ otherwise.
Similarly to the previous case we want  determine closed formulae and recursive relations for the coefficients that allow us to reduce the system. 
\begin{rmk}
 By Lemma \ref{simmetria11k}, it is enough to analyze the elements in the sets  $\Gamma_\nu^{i,j}$ for $i < n-i-j+3$ with the convention that a weight in $\Gamma_\nu^{i,j}$ gives a contribution to the coefficient equal to  $\Psi_i^j$ if $i > 0$ and equal to $F_0^j$ if $i=0$.
This simplification implies that $2i-2 < n-j+1=(e_j,\rho)$ and again without loss of generality we can suppose the coordinates of $ w\lambda -2ie_j$ are all non negative and the reduced form of $C_{w\lambda -2ie_j}$  can be again computed just by the action of the symmetric group on $w\lambda -2ie_j+\rho$.
\end{rmk}
\begin{rmk}\label{rho+peso1k}
Similarly to what observed in the case $\omega_{2k}$, the only possibility for a contribution to $\Lambda^{2|k, n}_{2|h}$ is by weights of the form $2e_1+e_2 + \dots + e_{2h+1} + \mu$ with $\mu \in \Omega_{0}^{k-h,n-2h-1}$, in particular $i=0$ and $w(e_1)=e_1$. By Lemma \ref{zerocontribution} we obtain:
\begin{equation}\label{coeff2h}
\Lambda^{2|k, n}_{2|h}=(-1)^{k-h}F_0^{1,n} \binom{n-k-h-1}{k-h}.
\end{equation}
\end{rmk}
Now we want find a recursive expansion of the coefficients $\Lambda^{2|k, n}_h$. Consider $\mu \in O_h^{2|k,n}:=\{w(\omega_1+\omega_{2k+1}-2je_1)|w(\omega_1+\omega_{2k+1}-2je_1)+\rho \sim \omega_{2h}+\rho \}$.
Let us suppose, first of all, that $h>1$. The only relevant cases are: 
\begin{itemize}
 \item\emph{Case 1: The fist two coordinates of $\mu$ are equal to 1.} Then  $\mu=(1,1,\mu')$ with $\mu' \in \Omega_{h-1}^{2|k-1,n-2}$
 \item \emph{Case 2; $\mu_1=1$ and $\mu_2=0$.}
This forces $\mu$ to be of the form $(1,0,2,1,\mu')$ with $\mu' \in \Omega_{h-2}^{k-1,n-4}$.
\item \emph{Case 3: $\mu_1=0$.} 
This forces $\mu = (0,2, \mu')$ with $\mu' \in \Omega_{h-1}^{k,n-2} $. 
Let us remark that this is equivalent to contract $(\mu_2, \dots , \mu_n)$ to a weight $w (\omega_1+ \omega_{2k+1})$ in $C_{n-1}$ conjugated to  $\omega_1+\omega_{2(h-1)+1}$. 
\end{itemize}
We can translate this enumerative analysis in the following recursive relation:
\[ \Lambda^{2|k, n}_h= \Lambda^{2 | k-1 , n-2}_{h-1}-\Lambda^{2| k , n-1}_{2| h-1} - (-1)^{k-h}F_0^{3,n} \binom{n-k-h-1}{k-h+1}.\]
If $h=1$, case 2 is not allowed, whereas we have to consider the additional case of weights of the form $(1,-1,2,\mu')$ where $\mu' \in \Omega_0^{k-1,n-3}$, leading us to the recurrence:
%
\[ \Lambda^{2|k, n}_1= \Lambda^{2| k-1 , n-2}_{0}-\Lambda^{2 |k , n-1}_{2| 0} +(-1)^{k}F_0^{3,n} \binom{n-k-2}{k-1}.\]
The zero case is again more difficult.
If we denote by $\Lambda^{2|k , n}_{0,i}$ the contribution to the coefficient $\Lambda^{2|k,n}_0$ for a fixed $i$, we have 
$\Lambda^{2|k,n}_0= \sum_{i=0}^{k+1} \Lambda^{2|k , n}_{0,i}.$
We want recover again a recurrence for $\Lambda_0^{2|k,n}$ summing up recurrences for each  $\Lambda^{2|k , n}_{0,i}$. Considering $\mu \in O_{0,i}^{2|k,n}:=\{w(\omega_1+\omega_{2k+1}-2je_1)|w(\omega_1+\omega_{2k+1}-2je_1) + \rho \sim \rho, \, j=i\}$ we have the following relevant cases:
\begin{itemize}
 \item \emph{Case 1: $\mu_1=0$} Then $\mu$ can be contracted to a weight $\mu' \in O_{0,i}^{2|k,n-1}$. 
 \item \emph{Case 2: $\mu_1=-1$ and the second to $\mu_2=1$.} Then $\mu=(-1,1,\mu')$ with $\mu' \in O_{0,i}^{2|k-1,n-2}$.
 \item  \emph{The first $j$ of $w \lambda$ coordinates are equal to $-1$.}. In this case $\mu + \rho$ is the vector $(n-1, \dots , n-2j, \dots )$. We must obtain a coordinate equal to $n$ because $\mu + \rho$ have to be conjugated to 0. This is possible only if $\mu=(-1,-1,2,\mu')$ (i.e. $i=0$ and $w e_1=e_3$) with $\mu' \in \Omega_{0}^{k-1,n-3}$.
 \item \emph{$\mu_1=2-2i$} This forces the second coordinate of $\mu$ to be equal to $1$, idem for the third and so on until the $2i-1$-th. Then $\mu=(2-2i, 1 , \dots, 1, \mu')$, where $\mu' \in \Omega_{0}^{k-i+1,n-2i}$.
\end{itemize}
We can translate what stated above in the following relations:
\[\Lambda_{0,i}^{2|k, \, n}= \Lambda_{0,i}^{2|k, \, n-1} - \Lambda_{0,i}^{2|k-1, \, n-2}+(-1)^{k-i+1} \Psi_{i}^{1,n} \binom{n-k-i}{k-i+1}  \qquad \mbox{if} \; \;  i>0, \]
 \[\Lambda_{0,0}^{2|k, \, n}= \Lambda_{0,0}^{2|k, \, n-1} - \Lambda_{0,0}^{2|k-1, \, n-2}+ (-1)^{k-1} F_0^{3,n} \binom{n-k-2}{k-1} \qquad \mbox{if} \; \;  i=0.\] 
Finally we can sum up all the contributions obtaining
\begin{equation}\label{Lambda2k0}\Lambda_{0}^{2|k, \, n}= \Lambda_{0}^{2|k, \, n-1} - \Lambda_{0}^{2|k-1, \, n-2}+ (-1)^{k-1}\binom{n-k-2}{k-1}F_0^{3, n}+ \sum_{i=1}^{k+1}(-1)^{k-i+1} \Psi_i^{1, n} \binom{n-k-i}{k-i+1}  \end{equation}
%
%
%
%
%
%
%
Imposing that the coefficients $\Lambda^{2|k,n}_{2|h}$ must cancel, an useful reduction to the triangular system can be found:
\begin{proposition}\label{magic}
 Let $\{R_i\}_{i \leq k}$ be the set of reduced recurrences. Then there exist coefficients $\Gamma^{k,n}_i$ and a family of integers $\{B_i^{k, n}\}_{i \leq k}$, with $B_k^{k,n}=1$ and $\sum_{i=h}^k(-1)^i \binom{n-i-1}{i}B_{i}^{k n }=0$, such that 
 \begin{equation}\label{n=2k+1}
 \sum B_i^{k,2k+1}R_i = \Lambda_{2|k}^{2|k, 2k+1} C_{2|k,2k+1}+ \Gamma_{k}^{k, 2k+1} C_{k,2k+1}  + \Gamma_0^{k, 2k+1} \left(C_{k-1, 2k+1}+ \dots + C_{0, 2k+1}\right),
\end{equation}
 and
 \begin{equation}\label{ngenerico}
\sum B_i^{k,n}R_i = \Lambda_{2|k}^{2|k, n} C_{2|k,n}+ \Gamma_{k+1}^{k, n} C_{k+1,n} + \Gamma_{k}^{k, n} C_{k,n}  + \Gamma_0^{k, n} \left(C_{k-1, n}+ \dots + C_{0, n}\right).
 \end{equation}
\end{proposition}
A proof of Proposition \ref{magic} can be found in Section 7. The fact that the coefficient of $C_{2|k}$ is exactly $\Lambda_{2|k}^{2|k,n}$ is a consequence of $B_k^{k,n}=1$.
In order to prove Reeder's Conjecture, we will find some explicit formulae for the coefficients $\Gamma_{0}^{k,n}$, $\Gamma_{k}^{k,n}$ and $\Gamma_{k+1}^{k,n}$.
In Section 7, Lemma \ref{FinediMondo} we prove that the following closed formula  for  $\Gamma_{0}^{k,n}$  holds:
\begin{equation*}
 \Gamma_0^{k,n}=\frac{(t^2+q)(t-q)(t^{2(n-2k)}-1)}{t^{n-2k+1}(t-1)}.
\end{equation*}
We can find a recursive expression for $\Gamma_k^{k,n}$ describing explicitly the contribution to the coefficient:
\begin{align*}
\Gamma_{k}^{k,n}=B_k^{k,n}\Lambda_k^{2|k,n}+B_{k-1}^{k,n}\Lambda_{k}^{2|k-1,n}=\\
\Lambda^{2| k-1 , n-2}_{k-1}-\Lambda^{2| k , n-1}_{2| k-1} + F_0^{3,n} \binom{n-2k-1}{1} + B_{k-1}^{k,n}\left(\Lambda^{2| k-2 , n-2}_{k-1}-\Lambda^{2| k-1 , n-1}_{2| k-1} -F_0^{3,n}\right)= \\
\Lambda^{2| k-1 , n-2}_{k-1}+B_{k-1}^{k,n}\Lambda^{2| k-2 , n-2}_{k-1} - \left[\Lambda^{2| k , n-1}_{2| k-1} + B_{k-1}^{k,n}\Lambda^{2| k-1 , n-1}_{2| k-1}\right]+ F_0^{3,n}\left[\binom{n-2k-1}{1} -B_{k-1}^{k,n} \right]=\\
\Lambda^{2| k-1 , n-2}_{k-1}+B_{k-1}^{k,n}\Lambda^{2| k-2 , n-2}_{k-1} - F_0^{1,n-1}\left[B_{k-1}^{k,n}-\binom{n-2k-1}{1}\right]+ F_0^{3,n}\left[\binom{n-2k-1}{1} -B_{k-1}^{k,n} \right]=\\
 =\\ \Gamma_{k-1}^{k-1,n-2} - \left(F_0^{2,n}+F_0^{3,n}\right)\left[B_{k-1}^{k,n}-\binom{n-2k-1}{1}\right]=\\
  =\\ \Gamma_{k-1}^{k-1,n-2} - \left(F_0^{2,n}+F_0^{3,n}\right)\left[B_{k-1}^{k,n}-\binom{n-2k}{1} +1\right]=\\
  \Gamma_{k-1}^{k-1,n-2} - \left[F_0^{2,n}+F_0^{3,n}\right].
\end{align*}
Moreover the only contribution to $\Gamma_{k+1}^{k,n}$ comes from $R_k$ and then 
\begin{equation*}
\Gamma_{k+1}^{k,n}=\Lambda_{k+1}^{2|k,n}=
\Lambda_{k}^{2|k-1,n-2}-\Lambda^{2| k , n-1}_{2| k} -F_0^{3,n} =  \Lambda^{2| k-1 , n-1}_{k} - \left[F_0^{2,n}+F_0^{3,n}\right] = \Gamma_{k}^{k-1,n-2} - \left[F_0^{2,n}+F_0^{3,n}\right].
\end{equation*}
%
Using an iterated process the following closed expression can be computed:
 \[\Gamma_{k}^{kn}=-\frac{(t-q)(1+qt^{2(n-2k)-1})}{t^{n-2k}}-\frac{(1-q^2t^{2n-2k-1})}{t^{n-1}} \frac{(t^{2k}-1)}{(t-1)},\]
\[\Gamma_{k+1}^{kn}=-\frac{(1-q^2t^{2(n-k-1)})}{t^{n-1}} \frac{(t^{2k+1}-1)}{(t-1)}.\]
\subsection{Reeder's Conjecture for $\lambda=\omega_1+\omega_{2k+1}$}
We already know that the following formula between the polynomials $C_{k, n}$ holds:
\[
 C_{k+1,n}=C_{k,n} \frac{q^4(q^{4(n-k+1)}-1)(q^{4(n-2k-1)}-1)(q^{4k-1}+1)}{(q^{4(k+1)}-1)(q^{4(n-2k+1)}-1)(q^{4(n-k)-1}+1)}. 
\]
We will denote this coefficient with $T_k^n$.
We recall that to prove the case $\omega_{2k}$ we shown the identity
$\Lambda_{k}^{k,n}C_{k}= - \Lambda_0^{1, n-2k-2} \left(C_{k-1}+ \dots + C_0\right)$.
%
%
%
%
Substituting it in relation (\ref{n=2k+1}), if $n=2k+1$ we obtain:
  \begin{align*}
  \Gamma_{2|k}^{k n} C_{2|k,n} = -\left[ \Gamma_{k}^{k, n} C_{k,n} + \Gamma_0^{k, n} \left(C_{k-1, n}+ \dots + C_{0, n}\right)\right] =
  -C_{k,n}\left[ \Gamma_{k}^{k, n} - \frac{\Gamma_0^{k, n} \Lambda_k^{k,n}}{\Lambda_0^{1,3}}  \right] 
 \end{align*}
 Set now 
\[D^n_k:=\Gamma_{k}^{k,n}-\Gamma_0^{k,n}=\frac{(q^2t^{2(k-1)}-1)(t^{2n-2k+1}-1)}{t^{n-1}(t-1)}\]
 Similarly to the previous case, substituting in relation (\ref{ngenerico}) we obtain:
 \begin{align*}
  \Gamma_{2|k}^{k, n} C_{2|k,n} =\\ -\left[ \Gamma_{k+1}^{k, n} C_{k+1,n} + \Gamma_{k}^{k, n} C_{k,n}  + \Gamma_0^{k, n} \left(C_{k-1, n}+ \dots + C_{0, n}\right)\right] = \\
  \\ -\left[ \Gamma_{k+1}^{k, n} C_{k+1,n} + \left(\Gamma_{0}^{k, n}+D_k^n\right) C_{k,n}  + \Gamma_0^{k, n} \left(C_{k-1, n}+ \dots + C_{0, n}\right)\right] = \\
 -\left[ \Gamma_{k+1}^{k, n} C_{k+1,n} + D_k^n C_{k,n}  + \Gamma_0^{k, n} \left(C_{k,n}+C_{k-1, n}+ \dots + C_{0, n}\right)\right] = \\
  -\left[ \left(\Gamma_{k+1}^{k n} - \frac{\Gamma_0^{k,n}\Lambda_{k+1}^{k+1,n}}{\Lambda_{0}^{1,n-2k}} \right) C_{k+1,n} + D_k^n C_{k,n} \right]= \\
   - C_{k+1,n} \left(\Gamma_{k+1}^{k, n} - \frac{\Gamma_0^{k,n}\Lambda_{k+1}^{k+1,n}}{\Lambda_{0}^{1,n-2}}  + \frac{D_k^n}{T_k^n}\right)
 \end{align*}
Comparing this result with formula (\ref{T2kk+1}), in both cases Reeder's Conjecture is now reduced to prove an univariate polynomial identity that we checked using the computer.
\begin{rmk}
 As in the case of the weights $\omega_{2k}$, Proposition \ref{magic} allow us to compute easily the closed formulae for the non specialized rational functions $C_{2|k}(q,t)$.
\end{rmk}

\section{Proof of Proposition 4.3}
 \begin{proposition}
\begin{equation}\label{rec}C_mc_m=\sum_{i=1}^{[\frac{m+1}{2}]}C_{m-2i+1}b_{i}+\sum_{i=1}^{[\frac{m}{2}]}C_{m-2i}b_{n-m+i+1}.
\end{equation}
\proof We will prove that the coefficient of $C_i,i<m$ in the right hand sides of \eqref{riscritta} and \eqref{rec} match.
\begin{lemma}\label{C0}
The coefficient of $C_0$ in the expression (\ref{riscritta}) is equal to the coefficient of $C_0$ in (\ref{rec}).
\proof
We want proceed by induction. Assume \eqref{rec} holds for $C_{m-h}$, $h\geq 1$.
Then 
\begin{equation}\label{cmmenoh}C_{m-h}c_{m-h}=\sum_{i=0}^{[\frac{m-h-1}{2}]}C_{m-h-2i-1}b_{i+1}+\sum_{i=1}^{[\frac{m-h}{2}]}C_{m-h-2i}b_{n-m+h+i+1}.\end{equation}
First consider the case $m=2k+1$.
If  $h$ is even, $m-h$ is odd and the coefficient of $C_0$ in \eqref{cmmenoh} is  $b_{(m-h-1)/2}$. If  $h$ is odd, $m-h$ is even the coefficient of $C_0$ in  \eqref{cmmenoh} is 
$b_{n+(m-h)/2+1}$.
Substituting into \eqref{riscritta}, the coefficient of $C_0$ is 
{\footnotesize
\begin{equation}\label{parziale}(-1)^{k}\Bigg (\binom{n-1+k}{k}b_1+\Gamma(2,n;k)\Bigg )+\sum_{s=0}^{k-1}(-1)^s\binom{n-m+s}{s}b_{n-k+s+1}+\sum_{r=1}^{k}(-1)^{r-1}\binom{n-m+r}{r}b_{k-r+1}.\end{equation}}
We want to show that this coefficients equals $b_{k+1}$, which is  the coefficient of $C_0$ in the right hand side of \eqref{rec}. This is in turn equivalent to prove the following equality.
\begin{equation}\label{eq}
\Gamma(2,n;k)=\sum_{s=0}^{k-1}(-1)^{s+k+1}\binom{n-2k-1+s}{s}b_{n-k+s+1}+\sum_{r=0}^{k-1}(-1)^{r+k}\binom{n-2k-1+r}{r}b_{k-r+1}.
\end{equation}
Set  $\Psi=(q+1)(q^{2n-2}+q^{2-2n})$; we remark that the following relation holds 
\begin{equation}\label{asd}
\Gamma(2,n;k)=\binom{n-k-2}{k-1}\Psi+\Gamma(2,n-2;k-1)+\Gamma(2,n-1;k).\end{equation}
Since $\Psi+b_{n-1}=b_n$, we have by induction
\begin{align*}
\Gamma(2,n;k)= \\
\binom{n-k-2}{k-1}\Psi+ \\ +\sum_{s=0}^{k-2}(-1)^{s+k}\binom{n-2k-1+s}{s}b_{n-k+s}+\sum_{r=0}^{k-2}(-1)^{r+k+1}\binom{n-2k-1+r}{r}b_{k-r}+\\
+\sum_{s=0}^{k-1}(-1)^{s+k+1}\binom{n-2k+s-2}{s}b_{n-k+s}+\sum_{r=0}^{k-1}(-1)^{r+k}\binom{n-2k+r-2}{r}b_{k-r+1}=\\
\binom{n-k-2}{k-1}\Psi+\binom{n-k-3}{k-1}b_{n-1}+\\ +\sum_{s=0}^{k-2}(-1)^{s+k}
\Bigg (\binom{n-2k-1+s}{s}-\binom{n-2k+s-2}{s}\Bigg )b_{n-k+s}\\+(-1)^kb_{k+1}+\sum_{r=1}^{k-1}(-1)^{r+k}\Bigg (\binom{n-2k+r-2}{r-1}+\binom{n-2k+r-2}{r}\Bigg)b_{k-r+1}=\\
\binom{n-k-2}{k-1}b_n-\binom{n-k-3}{k-2}b_{n-1}+\sum_{s=1}^{k-2}(-1)^{s+k}\binom{n-2k+s-2}{s-1}b_{n-k+s}+\\+\sum_{r=0}^{k-1}(-1)^{r+k}\binom{n-2k-1+r}{r}b_{k-r+1}=
\\\binom{n-k-2}{k-1}b_n-\binom{n-k-3}{k-2}b_{n-1}+\sum_{s=0}^{k-3}(-1)^{s+k+1}\binom{n-2k-1+s}{s}b_{n-k+s+1}+\\+\sum_{r=0}^{k-1}(-1)^{r+k}\binom{n-2k-1+r}{r}b_{k-r+1}=
\\\sum_{s=0}^{k-1}(-1)^{s+k+1}\binom{n-2k-1+s}{s}b_{n-k+s+1}+\sum_{r=0}^{k-1}(-1)^{r+k}\binom{n-2k-1+r}{r}b_{k-r+1}.
\end{align*}
Now assume  $m=2k$. 
Proceeding as above, the equality to prove is 
{\small
\begin{equation}\label{iop}b_{n-k+1}=
(-1)^{k-1}\Gamma(1,n;k) +\sum_{s=1}^{k-1}(-1)^{s-1}\binom{n-2k+s}{s}b_{n-k+s+1}+\sum_{r=0}^{k-1}(-1)^{r}\binom{n-2k+r}{r}b_{k-r}
\end{equation}}
or
$$\Gamma(1,n;k) =\sum_{s=0}^{k-1}(-1)^{k+s-1}\binom{n-2k+s}{s}b_{n-k+s+1}+\sum_{r=0}^{k-1}(-1)^{k+r}\binom{n-2k+r}{r}b_{k-r}
$$
As for \eqref{asd}, we have
\begin{equation}\label{asd2}\Gamma(1,n;k)=\binom{n-k-1}{k-1}\Psi+\Gamma(1,n-2;k-1)+\Gamma(1,n-1;k).\end{equation}
Again we have by induction
\begin{align*}
\Gamma(1,n;k)=\\
\binom{n-k-1}{k-1}\Psi+\sum_{s=0}^{k-2}(-1)^{s+k}\binom{n-2k+s}{s}b_{n-k+s}+\sum_{r=0}^{k-2}(-1)^{r+k+1}\binom{n-2k+r}{r}b_{k-r-1}\\
+
\sum_{s=0}^{k-1}(-1)^{s+k+1}\binom{n-2k+s-1}{s}b_{n-k+s}+\sum_{r=0}^{k-1}(-1)^{r+k}\binom{n-2k+r-1}{r}b_{k-r}=\\
\binom{n-k-1}{k-1}\Psi+\sum_{s=0}^{k-2}(-1)^{s+k}
\Bigg (\binom{n-2k+s}{s}-\binom{n-2k+s-1}{s}\Bigg )b_{n-k+s} +\binom{n-k-2}{k-1}b_{n-1}+ \\+\sum_{r=1}^{k-1}(-1)^{r+k}\Bigg (\binom{n-2k+r-1}{r}+\binom{n-2k+r-1}{r-1}\Bigg)b_{k-r}+(-1)^kb_k=\\
\binom{n-k-1}{k-1}\Psi+\sum_{s=0}^{k-2}(-1)^{s+k}
\binom{n-2k+s-1}{s-1}b_{n-k+s}\\+\binom{n-k-2}{k-1}b_{n-1}+\sum_{r=1}^{k-1}(-1)^{r+k}\binom{n-2k+r}{r}b_{k-r}+(-1)^kb_k=\\
\binom{n-k-1}{k-1}b_{n}-\binom{n-k-2}{k-2}b_{n-1}+\sum_{s=0}^{k-3}(-1)^{s+k+1}
\binom{n-2k+s}{s}b_{n-k+s+1}\\+\sum_{r=0}^{k-1}(-1)^{r+k}\binom{n-2k+r}{r}b_{k-r}=\\
=\sum_{s=0}^{k-1}(-1)^{k+s-1}\binom{n-2k+s}{s}b_{n-k+s+1}+\sum_{r=0}^{k-1}(-1)^{k+r}\binom{n-2k+r}{r}b_{k-r}.
\end{align*}
\endproof
\end{lemma}
Now we want prove the equality for the other coefficients. We start identifying such coefficients in equation (\ref{riscritta}) where we have substituted $C_{m-h}c_{m-h}$ with the corresponding expression (\ref{rec}).
\begin{align*}
 C_mc_m=\sum_{i=0}^{[\frac{m-1}{2}]}(-1)^{i}\binom{n-m+i}{i}C_{m-2i-1}c_{m-2i-1}+\sum_{i=1}^{[\frac{m}{2}]}(-1)^{i-1}\binom{n-m+i}{i}C_{m-2i}c_{m-2i}+ \\ +\sum_{i=0}^{[\frac{m-1}{2}]}(-1)^{i}C_{m-2i-1} \left(\binom{n-m+i}{i}b_1+\Gamma(2,n-m+2i+1;i)\right)+ \\ +\sum_{i=1}^{[\frac{m}{2}]}(-1)^{i-1}C_{m-2i}\Gamma(1,n-m+2i;i)=\\
 \sum_{i=0}^{[\frac{m-1}{2}]}(-1)^{i}\binom{n-m+i}{i}\left(\sum_{j=1}^{[\frac{m-2i}{2}]}C_{m-2i-2j}b_{j}+\sum_{j=1}^{[\frac{m-2i-1}{2}]}C_{m-2i-2j-1}b_{n-m+2i+j+2}\right)+ \\ +\sum_{i=1}^{[\frac{m}{2}]}(-1)^{i-1}\binom{n-m+i}{i}\left(\sum_{j=1}^{[\frac{m-2i+1}{2}]}C_{m-2i-2j+1}b_{j}+\sum_{j=1}^{[\frac{m-2i}{2}]}C_{m-2i-2j}b_{n-m+2i+j+1}\right)+ \\ +\sum_{i=0}^{[\frac{m-1}{2}]}(-1)^{i}C_{m-2i-1} \left(\binom{n-m+i}{i}b_1+\Gamma(2,n-m+2i+1;i)\right)+ \\ +\sum_{i=1}^{[\frac{m}{2}]}(-1)^{i-1}C_{m-2i}\Gamma(1,n-m+2i;i).
\end{align*}
Set now $m-h=2s$; then the coefficient of $C_h$ in the above expression is 
{\footnotesize
\begin{align*}
 \sum_{i+j=s}\sum_{i=0}^{ \lfloor \frac{m-1}{2}\rfloor} \sum_{j=1}^{ \lfloor \frac{m-2i}{2} \rfloor} (-1)^{i}\binom{n-m+i}{i} b_j  
 +\sum_{i+j=s} \sum_{i =1}^{\lfloor \frac{m-1}{2}\rfloor } \sum_{ j=1}^{ \lfloor \frac{m-2i}{2} \rfloor}(-1)^{i-1}\binom{n-m+i}{i}b_{n-h-j+1}  + (-1)^{s-1} \Gamma(1,n-h,s).
 \end{align*}}
If $m-h=2s+1$, the coefficient of $C_h$ is 
\begin{align*}
\sum_{i+j=s} \sum_{ i=0}^{ \lfloor \frac{m-1}{2}\rfloor} \sum_{ j=1}^{  \lfloor \frac{m-2i-1}{2} \rfloor} (-1)^{i}\binom{n-m+i}{i} b_{n-h-j+1} + \sum_{i+j=s} \sum_{i=1}^{ \lfloor \frac{m}{2}\rfloor} \sum_{ j =0}^{\lfloor \frac{m-2i-1}{2} \rfloor}(-1)^{i-1}\binom{n-m+i}{i}b_{j+1}+ \\  +(-1)^{s} \left(\binom{n-m+s}{s}b_1+\Gamma(2,n-h;s)\right).
 \end{align*}
We can rearrange the indices observing that the conditions on $j $ are redundant, obtaining
\begin{align*}
 \sum_{i=0}^{ s-1}  (-1)^{i}\binom{n-m+i}{i} b_{s-i} +  \sum_{i =1}^{s-1 } (-1)^{i-1}\binom{n-m+i}{i}b_{n-h-s+i+1} + (-1)^{s-1} \Gamma(1,n-h,s)
 \end{align*}
for $m-h=2s$ and
{\footnotesize
\begin{align*}
 \sum_{ i=0}^{ s-1} (-1)^{i}\binom{n-m+i}{i} b_{n-h-s+i+1}  +  \sum_{i=1}^{s} (-1)^{i-1}\binom{n-m+i}{i}b_{s-i+1} +(-1)^{s} \left(\binom{n-m+s}{s}b_1+\Gamma(2,n-h;s)\right)
 \end{align*}}
for $m-h=2s+1$. We will denote these expressions with  $\gamma^p(n,m,h)$ in the even case and with $\gamma^d(n,m,h)$ in the odd case.
For fixed $s$, $\gamma^p$ and $\gamma^d$ are translation invariants, i.e. $\gamma^{p/d}(n,m,h)=\gamma^{p/d}(n+1,m+1,h+1)$. 
Moreover, for fixed $s$, the same holds for the coefficients $b_{n-m+s+1}$ and $b_{s+1}$. 
We have to prove that $b_{n-m+s+1} = \gamma^p(n,k,h)$ and $b_{s+1} = \gamma^d(n,k,h)$. Using the translation invariance we can then reduce to the case $h=0$, that has been proved in Lemma \ref{C0}.
%
\endproof
\end{proposition}
\section{Proof of Propositions \ref{redrecursion} and \ref{magic}}
\begin{proposition}
 Let $R_i$ be the recurrence for $C_i(q,t)$ written in the reduced form. Then there exist a family of integers $\{A_i^{k, n}\}_{i \leq k}$ such that 
 \[\sum_{i=1}^k A_i^{k,n}R_i = \Lambda_k^{k, n} C_{k,n}(q, t)+ \Lambda_0^{1, n-2k+2} \left(C_{k-1, n}(q,t)+ \dots + C_{0, n}(q,t)\right)\]
\proof
Set
\begin{equation}
A_h^{k, n} =
\left\{
	\begin{array}{lll}
        0  & \mbox{if } h>k \mbox{ or } h\leq 0,\\
		1  & \mbox{if } h=k, \\
		A_{i-1}^{k-1,n-2}  & \mbox{if } k>h >1, \\
		\sum_{i=2}^k(-1)^i \binom{n-i-1}{i-1}A_{i-1}^{k-1,n-2} & \mbox{if } h=1.
	\end{array}
\right.
\end{equation}
We can rearrange the expression $\sum A_i^{k,n}R_i$ so that we are reduced to prove
$\sum _{i=h}^k A_i^{k, n } \Lambda_h^{i, n} = \Lambda_0^{1, n-2k + 2}$ for all $h<k$.
%
We start the proof of Proposition \ref{redrecursion} from the case $ h = 0$.  
Recalling that by definition $A_0^{k,n}=0$, we have to show that  
 \[\sum _{i=1}^k A_i^{k, n } \Lambda_0^{i, n} = \Lambda_0^{1, n-2k+ 2}.\]
Expanding $\sum _{i=1}^k A_i^{k, n } \Lambda_0^{i, n}$, using (\ref{L0}) and the definition of the $A_j^{k,n}$, we obtain:
 \begin{align*}
          \sum _{i=1}^k A_i^{k, n } \Lambda_0^{i, n} = \\
          \sum _{i=1}^k A_i^{k, n } \Lambda_0^{i, n-1} - \sum _{i=2}^k A_{i-1}^{k-1, n-2 } \Lambda_0^{i-1, n-2} +  \left(F_0^{2, n} + F_1^{1, n} \right) \left( \sum _{i=1}^k  (-1)^i A_i^{k, n }\binom{n-i-1}{i-1}\right) + \\ +  \sum _{i=2}^k  \left[ A_i^{k, n } \sum_{j=2}^{i}(-1)^{i-j+1} F_j^{1, n} \binom{n-i-j}{i-j} \right]= \\
           \sum _{i=1}^k A_i^{k, n } \Lambda_0^{i, n-1} - \sum _{t=1}^{k-1} A_{t}^{k-1, n-2 } \Lambda_0^{t, n-2} +  \left(F_0^{2, n} + F_1^{1, n} \right) \left[-A_1^{k, n} + \sum _{i=2}^k  (-1)^i A_i^{k, n }\binom{n-i-1}{i-1}\right] + \\ +  \sum _{j=2}^k F_j^{1, n}   \left[ \sum_{i=j}^k A_i^{k, n }(-1)^{i-j+1}\binom{n-i-j}{i-j} \right]=  \\
           \sum _{i=1}^k A_i^{k, n } \Lambda_0^{i, n-1} - \Lambda_0^{1, n-2k+2} +  \sum _{j=2}^k F_j^{1, n}   \left[ \sum_{t=1}^{k-j+1} A_{t+j-1}^{k, n }(-1)^{t}\binom{n-2j-t+1}{t-1} \right],
           \end{align*}
           where by inductive hypothesis we replaced $\sum _{t=1}^{k-1} A_{t}^{k-1, n-2 } \Lambda_0^{t, n-2}$ with $\Lambda_0^{1, n-2k+2}$. Now
           \begin{align*}
           \sum _{i=1}^k A_i^{k, n } \Lambda_0^{i, n-1} - \Lambda_0^{1, n-2k+2} +  \sum _{j=2}^k F_j^{1, n}   \left[ \sum_{t=1}^{k-j+1} A_{t+j-1}^{k, n }(-1)^{t}\binom{n-2j-t+1}{t-1} \right]=\\
             \sum _{i=1}^k A_i^{k, n } \Lambda_0^{i, n-1} - \Lambda_0^{1, n-2k+2} +  \sum _{j=2}^k F_j^{1, n}   \left[ \sum_{t=1}^{k-j+1} A_{t}^{k-j+1, n-2j+2}(-1)^{t}\binom{n-2j+2-t-1}{t-1} \right]=  \\
                \sum _{i=1}^k A_i^{k, n } \Lambda_0^{i, n-1} - \Lambda_0^{1, n-2k+2} - F_k^{1, n}.
           \end{align*}
           It can be shown by straightfoward computation that $A_h^{k, n}=A_h^{k, n-1}+A_h^{k-1, n-1}$, so that
\begin{align*}
         \sum _{i=1}^k A_i^{k, n } \Lambda_0^{i, n-1} - \Lambda_0^{1, n-2k+2} - F_k^{1, n}= 
         \sum _{i=1}^k \left[A_i^{k, n-1}+A_i^{k-1, n-1} \right]\Lambda_0^{i, n-1} - \Lambda_0^{1, n-2k+2} - F_k^{1, n}= \\
         \sum _{i=1}^k A_i^{k, n-1}\Lambda_0^{i, n-1}  + \sum _{i=1}^{k-1} A_i^{k-1, n-1}\Lambda_0^{i, n-1} - \Lambda_0^{1, n-2k+2} - F_k^{1, n}= \\ \Lambda_0^{1, n-2k+1}+ \Lambda_0^{1,n-2k+3} - \Lambda_0^{1, n-2k+2} - F_k^{1, n}= \\ \Lambda_0^{1, n-2k+1}- \sum_{j=2}^{n-2k+3} \left[F_0^{2,j}+F_1^{1,j} \right]+\sum_{j=2}^{n-2k+2} \left[F_0^{2,j}+F_1^{1,j} \right] - F_k^{1, n} = \\ \Lambda_0^{1, n-2k+1}- F_1^{1, n-2k+3}-F_0^{2,n-2k+3}-F_k^{1, n}.
           \end{align*}
Observing now that $ F_1^{1, n-2k+3}+F_0^{2,n-2k+3}+F_k^{1, n}= F_1^{1, n-2k+2}+F_0^{2,n-2k+2}$ we have:
          \begin{align*}
            \Lambda_0^{1, n-2k+1}- F_1^{1, n-2k+3}-F_0^{2,n-2k+3}-F_k^{1, n} =\\ \Lambda_0^{1, n-2k+1}  - \left[F_1^{1, n-2k+2}+F_0^{2,n-2k+2}\right] = \\ 
            -\sum_{j=2}^{n-2k+1} \left[F_0^{2,j}+F_1^{1,j} \right] -  \left[F_1^{1, n-2k+2}+F_0^{2,n-2k+2} \right] =\\ -\sum_{j=2}^{n-2k+2} \left[ F_1^{1, j}+F_0^{2,j} \right]= \Lambda_0^{1, n-2k+2}.
         \end{align*}
The case $k > h > 0$ can be deduced by some computations from the case $h=0$:
\begin{align*}
          \sum _{i=h}^k A_i^{k, n } \Lambda_h^{i, n}= 
           \sum _{i=1}^{k-h} A_{h+i}^{k, n } \Lambda_h^{h+i n}+A_h^{k, n}\Lambda_h^{h, n} =\\
           \sum _{i=1}^{k-h} A_{h+i}^{k, n }\left[ \Lambda_0^{i,n-2h} + (-1)^i \binom{n-2h-i}{i}\Lambda_h^{h, n}\right]+A_h^{k, n}\Lambda_h^{h, n} = \\
           \sum _{i=1}^{k-h} A_{i}^{k-h, n-2h }\Lambda_0^{i, n-2h} + \Lambda_h^{h, n} \left[ A_h^{k, n} + \sum_{i=1}^{k-h}(-1)^i \binom{n-2h-i}{i}A_{h+i}^{k, n}\right] = \\
           \Lambda_0^{1,  n-2k+2} + \Lambda_h^{h, n} \left[ A_{h}^{k, n} + \sum_{i=1}^{k-h}(-1)^i \binom{n-2h-i}{i}A_{i}^{k-h,n-2h}\right]=\\ 
             \Lambda_0^{1,  n-2k+2} + \Lambda_h^{h, n} \left[ A_1^{k-h+1, n-2h+2} -\sum_{t=2}^{k-h+1}(-1)^{t} \binom{n-2h+2-t-1}{t-1}A_{t-1}^{k-h, n-2h}\right]= \\
           \Lambda_0^{1,  n-2k+2}.
         \end{align*}
         \endproof
   \end{proposition}
The proof of Proposition 7.4 is very similar to the previous one.
\begin{proposition}
 Let $\{R_i\}_{i \leq k}$ be the set of reduced recurrences, then there exists a family of integers $\{B_i^{k n}\}_{i \leq k}$ such that 
 \begin{align*}
  \sum B_i^{k,2k+1}R_i = \Gamma_{2|k}^{k, 2k+1} C_{2|k,2k+1}+ \Gamma_{k}^{k, 2k+1} C_{k,2k+1}  + \Gamma_0^{k, 2k+1} \left(C_{k-1, 2k+1}+ \dots + C_{0, 2k+1}\right).
\end{align*}
 and 
 \begin{align*} \sum B_i^{k,n}R_i = \Gamma_{2|k}^{k, n} C_{2|k,n}+ \Gamma_{k+1}^{k, n} C_{k+1,n} + \Gamma_{k}^{k, n} C_{k,n}  + \Gamma_0^{k, n} \left(C_{k-1, n}+ \dots + C_{0, n}\right).
 \end{align*}
\proof
Set
\begin{equation}
B_h^{k, n} =
\left\{
	\begin{array}{lll}
        0  & \mbox{if } h>k \mbox{ or } h<0, \\
		1  & \mbox{if } h=k, \\
		B_{i-1}^{k-1, n-2}  & \mbox{if } k>h >0, \\
		\sum_{i=1}^k(-1)^i \binom{n-i-1}{i}B_{i}^{k, n } & \mbox{if } h=0.
	\end{array}
\right.
\end{equation}
%
%
%
%
Using the definition, some straightforward computations show that
\begin{lemma}
 \begin{enumerate}
  \item 
  \begin{equation} \label{idBk}
 \sum_{i=h}^k(-1)^{i-h} \binom{n-i-h-1}{i-h}B_{i}^{k, n }=0,
 \end{equation}
\item \[B_h^{k+1, n+1}=B_h^{k+1, n}+B_h^{k, n},\]
  \item  
\begin{equation}\label{tridBk}
 \sum_{i=h}^k(-1)^{i-h} \binom{n-i-h-2}{i-h}B_{i}^{k,n} =0 \qquad (h<k-1),
\end{equation}
\item   $B_j^{k,2k}=0$ for all $j<k$.
\end{enumerate}
\end{lemma}
We will denote by $\Gamma_{2|h}^{k,n}$ and by $\Gamma_{h}^{k,n}$ respectively the coefficients of $C_{2|h}$ and of $C_h$ in $\sum_{i=1}^n B_i^{k,n}R_i$.
As a consequence of Lemma \ref{idBk}, substituting the formula (\ref{coeff2h}) and using 1),  we have 
for all $0 \leq h < k$ 
\[\Gamma^{k,n}_{2|h}=\sum_{i=h}^k B_i^{k, n} \Lambda_{2|h}^{2|i, n} =0.\]
In order to prove Proposition \ref{magic} it is enough now to show  $\Gamma^{k,n}_{h} =  \Gamma^{k,n}_{0}$ for all $h<k$.
More precisely, we want prove an iterative formula to compute the coefficients $\Gamma_{h}^{k,n}$. 
Let us start from the case $k-1\geq h>1$. We have:
\begin{align*}
 \Gamma_{h}^{k,n} = \sum _{i=h-1}^k B_i^{k, n } \Lambda_h^{2|i, n} = \\
 \sum_{i=h-1}^kB_i^{k,n} \Lambda_{h-1}^{2|i-1, n-2} - \sum_{i=h-1}^kB_i^{k,n} \Lambda_{2|h-1}^{2|i, n-1}- \left[\sum_{i=h-1}^k (-1)^{i-h+1}B_i^{k,n} \binom{n-i-h-1}{i-h+1}\right]F_0^{3,n}=\\
 \sum_{i=h-1}^kB_{i-1}^{k-1,n-2} \Lambda_{h-1}^{2|i-1, n-2} - \sum_{i=h-1}^kB_i^{k,n} \left[(-1)^{i-h+1}\binom{n-i-h-1}{i-h+1}\right] F_0^{1,n-1}+ \\ -\left[\sum_{i=h-1}^k (-1)^{i-h+1}B_i^{k,n} \binom{n-i-h-1}{i-h+1}\right]F_0^{3,n}=\\
  \sum_{t=h-2}^{k-1}B_{t}^{k-1,n-2} \Lambda_{h-1}^{2|t, n-2}=\Gamma^{k-1,n-2}_{h-1}.
 \end{align*}
The result is exactly the same if $n=2k+1$, but in this case the much simpler recursions lead to simpler computations.
 Similar computations in the case $h=1$ lead to $\Gamma_1^{k,n}=\Gamma_0^{k-1,n-2}$.
If $n=2k+1$ we need to expand the $B_i^{k,n}$ using Lemma \ref{idBk}: 
\begin{align*}
 \Gamma_{1}^{k,n} = \sum _{i=0}^k B_i^{k, n } \Lambda_1^{2|i, n} = \\
 \sum_{i=1}^kB_i^{k,n} \Lambda_{1}^{2|i-1, n-2} - \sum_{i=0}^{k-1}B_i^{k,n} \Lambda_{2|0}^{2|i, n-1}+ \left[\sum_{i=0}^k (-1)^{i-2}B_i^{k,n} \binom{n-i-2}{i-1}\right]F_0^{3,n}=\\
 \sum_{t=0}^{k-1}B_{t+1}^{k,n} \Lambda_{0}^{2|t, n-2} - \sum_{i=0}^{k-1}\left[B_{i}^{k-1,n-1}+B_{i}^{k,n-1}\right] \Lambda_{2|0}^{2|i, n-1}=\\
  \sum_{t=0}^{k-1}B_{t}^{k-1,n-2} \Lambda_{0}^{2|t n-2}-\sum_{i=0}^{k-1}B_{i}^{k-1,n-1} \Lambda_{2|0}^{2|i, n-1} = \\ \Gamma^{k-1,n-2}_{0} -\Gamma^{k-1,n-2}_{2|0} = \Gamma^{k-1,n-2}_{0}
 \end{align*}
We have finally only to compute the coefficient $\Gamma_0^{k,n}$. Recalling equation (\ref{Lambda2k0})
 and using such expression to expand $\sum _{i=1}^k B_i^{k, n } \Lambda_0^{i, n}$, we obtain:
\begin{align*}
 \Gamma_0^{k,n}=\sum _{i=0}^k B_i^{k, n } \Lambda_0^{2|i, n} = \\
 \sum _{i=0}^k B_i^{k, n } \Lambda_0^{2|i, n-1} - \sum _{i=0}^{k-1} B_{i+1}^{k, n  } \Lambda_0^{2|i, n-2} +  \sum _{i=1}^k  (-1)^{i-1} B_i^{k, n }\binom{n-i-2}{i-1} F_0^{3, n}  + \\ +  \sum _{i=0}^k  \left[ B_i^{k, n } \sum_{j=1}^{i+1}(-1)^{i-j+1} \Psi_j^{1, n} \binom{n-i-j}{i-j+1} \right] = \\
\sum _{i=0}^k B_i^{k, n } \Lambda_0^{2|i, n-1} - \sum _{i=0}^{k-1} B_{i+1}^{k, n  } \Lambda_0^{2|i, n-2} +  F_0^{3, n}  \left( \sum _{i=1}^k  (-1)^{i-1} B_i^{k, n }\binom{n-i-2}{i-1}\right) + \\ +  \sum _{j=1}^{k+1} \Psi_j^{1,n} \left(  \sum_{i=j-1}^{k}(-1)^{i-j+1}\binom{n-i-j}{i-j+1}B_i^{k, n } \right)= \\
\sum _{i=0}^k B_i^{k, n } \Lambda_0^{2|i, n-1} - \sum _{i=0}^{k-1} B_{i+1}^{k, n  } \Lambda_0^{2|i, n-2} +  F_0^{3, n}  \left( \sum _{i=1}^k  (-1)^{i-1} B_i^{k, n }\binom{n-i-2}{i-1}\right) + \\ +  \sum _{j=1}^{k+1} \Psi_j^{1,n} \left(  \sum_{t=0}^{k-j+1}(-1)^{t}\binom{n-2j-t+1}{t}B_{t+j-1}^{k, n } \right)= \\
\sum _{i=0}^k B_i^{k, n } \Lambda_0^{2|i, n-1} - \sum _{i=0}^{k-1} B_{i+1}^{k, n  } \Lambda_0^{2|i, n-2} +  F_0^{3, n}  \left( \sum _{s=0}^{k-1}  (-1)^{s} B_{s}^{k-1, n-2 }\binom{(n-2) -s-1}{s}\right) + \\ +  \sum _{j=1}^{k} \Psi_j^{1,n} \left(  \sum_{t=0}^{k-j+1}(-1)^{t}\binom{n-2j-t+1}{t}B_{t}^{k-j+1, n-2j+2 } \right)+\Psi_{k+1}^{1,n}= \\
\sum _{i=0}^k B_i^{k, n } \Lambda_0^{2|i, n-1} - \sum _{i=0}^{k-1} B_{i+1}^{k, n  } \Lambda_0^{2|i, n-2} +\Psi_{k+1}^{1,n}.
  \end{align*}
 Now expanding $B_i^{k,n}$ by Lemma \ref{idBk} and using the properties of the $B_i^{k,n}$ we obtain
\begin{align*}
 \sum _{i=0}^k B_i^{k, n } \Lambda_0^{2|i, n-1} - \sum _{i=0}^{k-1} B_{i+1}^{k, n  } \Lambda_0^{2|i, n-2} +\Psi_{k+1}^{1,n}=\\
 \sum _{i=0}^k \left[B_i^{k, n-1 } + B_i^{k-1,n-1}\right] \Lambda_0^{2|i, n-1} - \sum _{i=0}^{k-1} B_{i}^{k-1, n -2 } \Lambda_0^{2|i, n-2} +\Psi_{k+1}^{1,n} =\\
 \sum _{i=0}^k B_i^{k, n-1 } \Lambda_0^{2|i, n-1} + \sum _{i=0}^{k-1}  B_i^{k-1n-1} \Lambda_0^{2|i, n-1} - \sum _{i=0}^{k-1} B_{i}^{k-1, n -2 } \Lambda_0^{2|i, n-2} +\Psi_{k+1}^{1,n} =\\
 \Gamma_0^{k,n-1}+\left(\Gamma_0^{k-1,n-1}-\Gamma_0^{k-1,n-2}\right)+\Psi_{k+1}^{1,n}.
  \end{align*}
The case $n=2k+1$ is very similar and leads to
  $\Gamma_0^{k,n} = \Gamma_0^{k-1,n-1}-\Gamma_0^{k-1,n-2}+\Psi_{k+1}^{1,n}.$
\begin{lemma}\label{FinediMondo}
Set 
$\Psi(n,k)=\sum_{i=1}^{k+1}\Psi_i^{k-i+2, n} + F_0^{k+2,n}.$
 The following identities hold:
\begin{enumerate}
  \item \[\Gamma_0^{k,n}=\sum_{j=2k+1}^{n}\Psi(j,k).\]
  \item \[\Psi(n,k)=\frac{(t^2+q)(t-q)(t^{2(n-2k)-1}+1)}{t^{n-2k+1}}.\]
  \item \[\Gamma_0^{k,n}=\Gamma_0^{k-1,n-2}=\Gamma_1^{k,n}.\]
  \item \[\Gamma_0^{k,n}=\frac{(t^2+q)(t-q)(t^{2(n-2k)}-1)}{t^{n-2k+1}(t-1)}.\]
 \end{enumerate}
\proof
We will apply an inductive reasoning. Observe that, by definition
$\Psi(n,k)=\Psi(n-1,k-1)+\Psi_{k+1}^{1,n}$. We remark that Lemma can be easily checked for $k=1$ and $n=3$:\[\Gamma_0^{1,3}= F_0^{3,3}+\Psi_1^{2,3}+\Psi_2^{1,3}=\frac{(t^2+q)(t-q)(t+1)}{t^2}.\]
\begin{enumerate}
\item Suppose there exists a minimal pair $(n,k)$ such that $1)$ does not hold. If $n=2k+1$ 
 \[\Gamma_0^{k,n}=\Gamma_0^{k-1,n-1}-\Gamma_0^{k-1,n-2}+\Psi_{k+1}^{1,n}=\Psi(n-1,k-1)+\Psi_{k+1}^{1,n}=\Psi(n,k).\]
If $n$ is generic, such a pair cannot exists by a very similar argument:
 \begin{align*}
  \Gamma_0^{k,n}=\Gamma_0^{k-1,n-1} + (\Gamma_0^{k-1,n-1}-\Gamma_0^{k-1,n-2})+\Psi_{k+1}^{1,n}=\\\Gamma_0^{kn-1}+\Psi(n-1,k-1)+\Psi_{k+1}^{1,n}=\Gamma_0^{k,n-1}+\Psi(n,k)= \sum_{j=2k+1}^n \Psi(j,k).
 \end{align*}
\item  
By induction we obtain:
\begin{align*}
 \Psi(n,k)=\Psi(n-1,k-1)+\Psi_{k+1}^{1,n}=\\
 \frac{(t^2+q)(t-q)(t^{2(n-2k+1)-1}+1)}{t^{n-2k+2}} +\frac{(t^2 + q)(q - t)(t^{2n}-t^{4k})(t - 1)}{t^{n+2k+2}}=\\
 \frac{(t^2+q)(t-q)(t^{2(n-2k)-1}+1)}{t^{n-2k+1}}
\end{align*}
\item First of all let us observe that, by $2)$, we have $\Psi(n,k)=\Psi(n-2,k-1)$. Now, if $n=2k+1$, the result is an immediate consequence of $1)$, because coefficients $\Gamma_0^{k,2k+1}$ and $\Gamma_0^{k-1,n-1}$ are equal to $\Psi(2k+1,k)$ and $\Psi(2(k-1)+1,k-1)$ respectively.
Otherwise by induction we have
 \begin{align*}
  \Gamma_0^{k,n}=\Gamma_0^{k-1,n-1} + (\Gamma_0^{k-1,n-1}-\Gamma_0^{k-1,n-2})+\Psi_{k+1}^{1,n}=\Gamma_0^{k,n-1}+\Psi(n,k)=\Gamma_0^{k-1,n-3}+\Psi(n,k).
 \end{align*}
Now we have the following equivalences 
 \[\Gamma_0^{k-1,n-3}+\Psi(n,k) = \Gamma_0^{k-1,n-2} \Longleftrightarrow \Psi(n,k) = \Gamma_0^{k-1,n-2}- \Gamma_0^{k-1,n-3} = \Psi(n-2,k-1)\]
 and the statement comes again from $\Psi(n,k)=\Psi(n-2,k-1)$.
\item By $3)$ it is enough prove the statement for $k=1$.
Now by induction the statement follows immediately from 1).
\end{enumerate}
\endproof
\end{lemma}
As an immediate Corollary, recalling that $\Gamma_h^{k,n}=\Gamma_{h-1}^{k-1, n-2}$ if $h>0$ and reasoning by induction we obtain that 
\[\Gamma_0^{k,n}=\Gamma_1^{k,n}=\Gamma_0^{k-1,n-2}=\Gamma_j^{k-1,n-2}=\Gamma_{j+1}^{k,n}\]
for $k-1 > j \geq 0$ and this proves the statement of Proposition \ref{magic}. 
\endproof
\end{proposition}

\end{document}